\documentclass[11pt]{article}
\RequirePackage[OT1]{fontenc}
\RequirePackage{amsthm,amsmath}
\RequirePackage{hyperref}
\usepackage{amsfonts}
\usepackage[mathscr]{eucal}
\usepackage{amssymb}
\usepackage{amstext}
\usepackage{parskip}
\usepackage{fullpage}
\usepackage{graphicx}
\usepackage{subcaption}

\RequirePackage[numbers]{natbib}
\usepackage{amsfonts,amssymb,amsthm,amsmath,enumerate}

{\end{eqnarray}\end{subequations}\hskip-4.0pt}

\def\fatnorm#1{|\kern-.2ex|\kern-.2ex| #1 |\kern-.2ex|\kern-.2ex|}
\newcommand{\twonorm}[1]{\left\lVert#1\right\rVert_2}

\newcommand{\fnorm}[1]{\left\lVert#1\right\rVert_{F}}
\newcommand{\norm}[1]{\left\lVert#1\right\rVert}
\newcommand{\abs}[1]{\left\lvert#1\right\rvert}

\newcommand{\Sp}{\mathbb{S}}
\newcommand{\T}{\mathcal{T}}

\newcommand{\cov}{\textsf{Cov}}

\newcommand{\off}{{\rm off}}

\newcommand{\half}{\ensuremath{\frac{1}{2}}}
\newcommand{\inv}[1]{\frac{1}{#1}}
\def\conv{\mathop{\text{\rm conv}\kern.2ex}}

\newcommand{\ip}[1]{\;\langle{\,#1\,}\rangle\;}

\newcommand{\onenorm}[1]{\ensuremath{\left|#1\right|_1}}
\newcommand{\offone}[1]{\ensuremath{\left|#1\right|_{1,\off}}}

\newcommand{\expct}[1]{\ensuremath{\mathbb E}\left(#1\right)}
\newcommand{\silent}[1]{}
\newcommand{\mvec}[1]{\rm{vec}\{\,#1\,\}}

\newcommand{\ve}{\varepsilon}

\newcommand{\ksump}{\mathcal{K}_{\mathbf{p}}}

\newcommand{\drp}{\mathcal{D}}

\def\qed{\hskip1pt $\;\;\scriptstyle\Box$}
\def\Ber{\mathop{\text{Bernoulli}\kern.2ex}}
\def\supp{\mathop{\text{supp}\kern.2ex}}
\def\corr{\mathop{\text{corr}\kern.2ex}}
\def\prec{\mathop{\text{precision}\kern.2ex}}
\def\recall{\mathop{\text{recall}\kern.2ex}}
\def\cov{\mathop{\text{Cov}\kern.2ex}}
\def\mnorm{\mathcal{N}_{f,m}\kern.2ex}
\def\var{\mathop{\text{Var}\kern.2ex}}
\def\ess{\mathop{\text{ess}\kern.2ex}}
\def\dom{\mathop{\text{dom}\kern.2ex}}
\def\lin{\mathop{\text{lin}\kern.2ex}}

\newcommand{\func}[1]{\ensuremath{\mathrm{#1}}}
\newcommand{\diag}{\func{diag}}
\newcommand{\offd}{\func{offd}}

\let\hat\widehat
\let\tilde\widetilde

\newcommand{\tr}{{\rm tr}}

\def\E{{\mathbb E}}

\def\MG{{\mathcal G}}

\def\supp{\mathop{\text{\rm supp}\kern.2ex}}
\def\argmin{\mathop{\text{arg\,min}\kern.2ex}}

\DeclareMathOperator*{\Span}{span}

\newcommand{\prob}[1]{\ensuremath{\mathbb P}\left(#1\right)}

\newcommand{\beq}{\begin{equation}}
\newcommand{\eeq}{\end{equation}}
\newcommand{\ben}{\begin{eqnarray}}
\newcommand{\een}{\end{eqnarray}}
\newcommand{\bnum}{\begin{enumerate}}
\newcommand{\enum}{\end{enumerate}}
\newcommand{\bit}{\begin{itemize}}
\newcommand{\eit}{\end{itemize}}
\newcommand{\bens}{\begin{eqnarray*}}
\newcommand{\eens}{\end{eqnarray*}}
\newcommand{\X}{\boldsymbol{\mathscr{X}}}
\newcommand{\bA}{\boldsymbol{\mathscr{A}}}
\newcommand{\Y}{\boldsymbol{\mathscr{Y}}}
\newcommand{\bX}{\ensuremath{{\bf {X}}}}
\newcommand{\bY}{\ensuremath{{\bf {Y}}}}
\newcommand{\R}{{\mathbb R}}
\newcommand{\bp}{\mathbf{p}}

\newcommand{\Sc}{\ensuremath{S^c}}

\newcommand{\Ball}{{B}}

\def\qed{\hskip1pt $\;\;\scriptstyle\Box$}

\newenvironment{proofof2}{\hskip10pt}{\qed\vskip5pt}

\def\fatnorm#1{|\kern-.2ex|\kern-.2ex| #1 |\kern-.2ex|\kern-.2ex|}
\def\conv{\mathop{\text{\rm conv}\kern.2ex}}
\def\qed{\hskip1pt $\;\;\scriptstyle\Box$}
\def\Ber{\mathop{\text{Bernoulli}\kern.2ex}}
\def\supp{\mathop{\text{supp}\kern.2ex}}
\def\corr{\mathop{\text{corr}\kern.2ex}}
\def\prec{\mathop{\text{precision}\kern.2ex}}
\def\recall{\mathop{\text{recall}\kern.2ex}}
\def\cov{\mathop{\text{Cov}\kern.2ex}}
\def\mnorm{\mathcal{N}_{f,m}\kern.2ex}
\def\var{\mathop{\text{Var}\kern.2ex}}
\def\ess{\mathop{\text{ess}\kern.2ex}}
\def\dom{\mathop{\text{dom}\kern.2ex}}
\def\lin{\mathop{\text{lin}\kern.2ex}}

\let\hat\widehat
\let\tilde\widetilde

\def\E{{\mathbb E}}

\def\supp{\mathop{\text{\rm supp}\kern.2ex}}
\def\argmin{\mathop{\text{arg\,min}\kern.2ex}}

\graphicspath{FiguresTalk/}
 \graphicspath{{Figures/}}

\newtheorem{theorem}{Theorem}[section]

\newtheorem{lemma}[theorem]{Lemma}
\newtheorem{proposition}[theorem]{Proposition}

\newtheorem{definition}[theorem]{Definition}

\newtheorem{corollary}[theorem]{Corollary}

\def\qed{\hskip1pt $\;\;\scriptstyle\Box$}

\begin{document}

\title{Finite sample rates of convergence for the Bigraphical and Tensor graphical Lasso estimators}

\author{Shuheng Zhou \ \ \ \ \ \ \ \ \ \ \  \ \ \ \ \ \ \ \ \ \ \ \ \ \ Kristjan Greenewald  \\
University of California, Riverside \ \ \ \ \ \ \ MIT-IBM Watson AI Lab, Cambridge, MA}

\date{}

\maketitle

\begin{abstract}
Many modern datasets exhibit dependencies among observations as well as
variables. A decade ago, Kalaitzis et. al. (2013) proposed the Bigraphical
Lasso, an estimator for precision matrices of matrix-normals based on the
Cartesian product of graphs; they observed that the associativity of the
Kronecker sum yields an approach to the modeling of datasets organized into 3
or higher-order tensors. Subsequently, Greenewald, Zhou and Hero (2019)
explored this possibility to a great extent, by introducing the tensor
graphical Lasso (TeraLasso) for estimating sparse $L$-way decomposable
inverse covariance matrices for all $L \ge 2$, and showing the rates of
convergence in the Frobenius and operator norms for estimating this class of
inverse covariance matrices for sub-gaussian tensor-valued data.  In this
paper, we provide sharper rates of convergence for both Bigraphical and
TeraLasso estimators for inverse covariance matrices.  This improves upon the
rates presented in GZH 2019.  In particular, (a) we strengthen the bounds for
the relative errors in the operator and Frobenius norm by a factor of
approximately $\log p$; (b) Crucially, this improvement allows for finite
sample estimation errors in both norms to be derived for the two-way
Kronecker sum model.  This closes the gap between the low single-sample error
for the two-way model as observed in GZH 2019 and the lack of theoretical
guarantee for this particular case.  The two-way regime is important because
it is the setting that is the most theoretically challenging, and
simultaneously the most common in applications.  Part of this work was
presented as a short conference paper in IEEE International Symposium on
Information Theory (ISIT 2024).  In the current paper, we elaborate on the
Kronecker Sum model, highlight the proof strategy and provide full proofs of
all main theorems.  Normality is not needed in our proofs; instead, we
consider subgaussian ensembles and derive tight concentration of measure
bounds, using tensor unfolding techniques.
\end{abstract}

\section{Introduction}
\label{Sec:Not}
Matrix and tensor-valued data with complex dependencies are ubiquitous in modern statistics and machine learning, flowing from sources as diverse as medical and radar imaging modalities, spatial-temporal and meteorological data collected from sensor networks and weather stations, and biological, neuroscience and spatial gene expression data aggregated over trials and time points. Learning useful structures from these large scale, complex and high-dimensional data in the low sample regime is an important task.
Undirected graphs are often used to describe high dimensional distributions.
Under sparsity conditions, the graph can be estimated using $\ell_1$-penalization methods, such as the graphical Lasso (GLasso) \citep{FHT07} and multiple nodewise regressions \citep{MB06}.
Under suitable conditions, including independence among samples, such approaches yield consistent and sparse estimation in terms of graphical structure and fast convergence rates with respect to the operator and Frobenius norm for the covariance matrix and its inverse. The independence assumptions substantially simplify mathematical derivations but tend to be very restrictive.

\begin{figure}
\vspace{-0.3in}
\begin{center}
  \begin{tabular}{c}
\vspace{-0.4in}
\includegraphics[width=0.95\textwidth,angle=0]{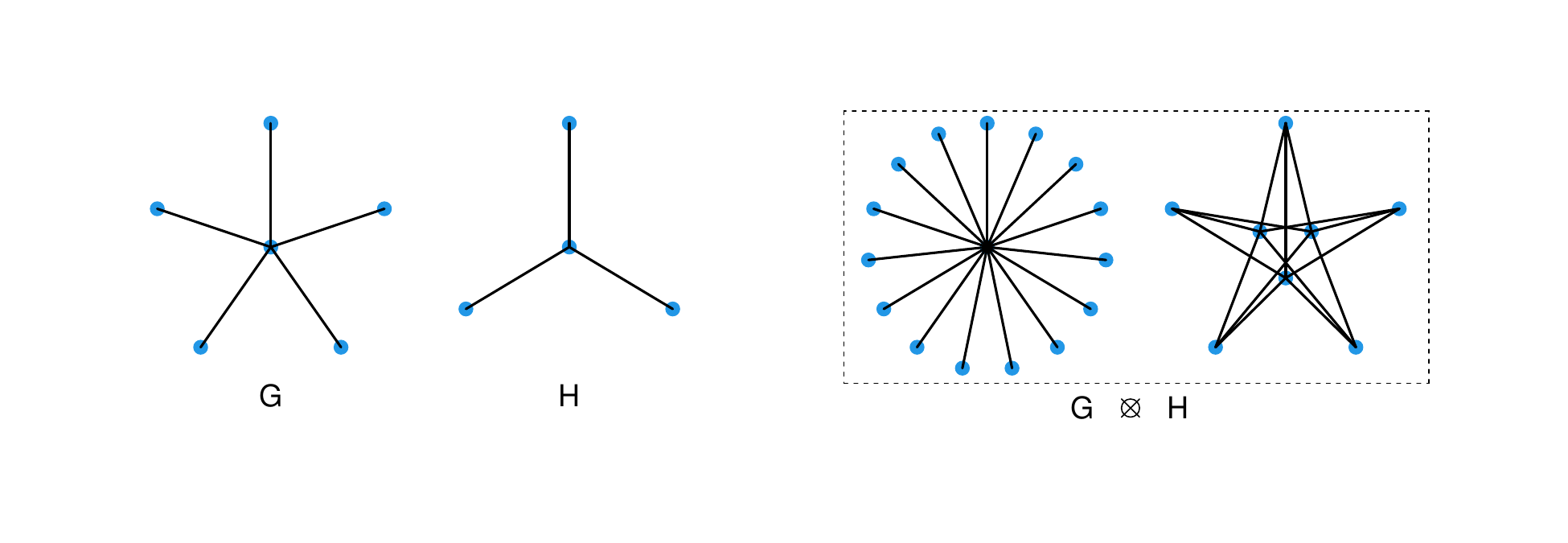}
\end{tabular}
\end{center}
  \caption{
  The Kronecker product of two graphs $G$ (corresponding to $A^{-1}$)
  and $H$ (corresponding to $B^{-1}$) is
  the graph whose adjacency matrix is the tensor product of the
  adjacency matrices of  $G$ and $H$~\citep{Weich62}.
  Observation: Estimating their Kronecker product directly following the classical
$p$-variate Gaussian graphical modeling approach will be costly in terms
of both computation and the sample requirements.}
\label{fig::prodgraph}
\vspace{-0.1in}
\end{figure}

To remedy this, recent work has demonstrated another regime where
further improvements in the sample size lower bounds are possible
under additional structural assumptions, which arise naturally in the
above mentioned contexts for data with complex dependencies. For example, the matrix-normal model~\citep{Dawid81} as studied in \cite{AT10},\cite{LT12},\cite{THZ13} and
\cite{Zhou14a} restricts the topology of the graph to tensor product
graphs where the precision matrix $A^{-1} \otimes B^{-1}$ corresponds
to a Kronecker product over two component graphs (cf. Figure~\ref{fig::prodgraph}). In~\cite{Zhou14a}, the author showed that one can estimate the
covariance and inverse covariance matrices well using only one
instance from the matrix variate normal distribution. However, such a
normality assumption is also not needed, as elaborated in a recent
paper by the same author in~\cite{Zhou24}.
More specifically, while the
precision matrix encodes conditional independence relations for
Gaussian distributions, for the more general sub-gaussian matrix
variate model, this no longer holds. However, the inverse covariance
matrix still encodes certain zero correlation relations between the
residual errors and the covariates in a regression model, analogous to
the Gaussian graphical models~\citep{laur96}. See~\cite{Zhou24}, where
such regression model is introduced for sub-gaussian matrix variate
data. See also~\cite{Horns19},~\cite{GPZ17},~\cite{fan19a}, and references therein for
recent applications of matrix variate models in genomics, neuroimaging
and political science.

\begin{figure}
  \begin{center}
          \vspace*{-1.5in}
          \vspace{-0.3in}
                    \hspace*{-0.2in}
\begin{tabular}{c}
    \includegraphics[width=4.5in]{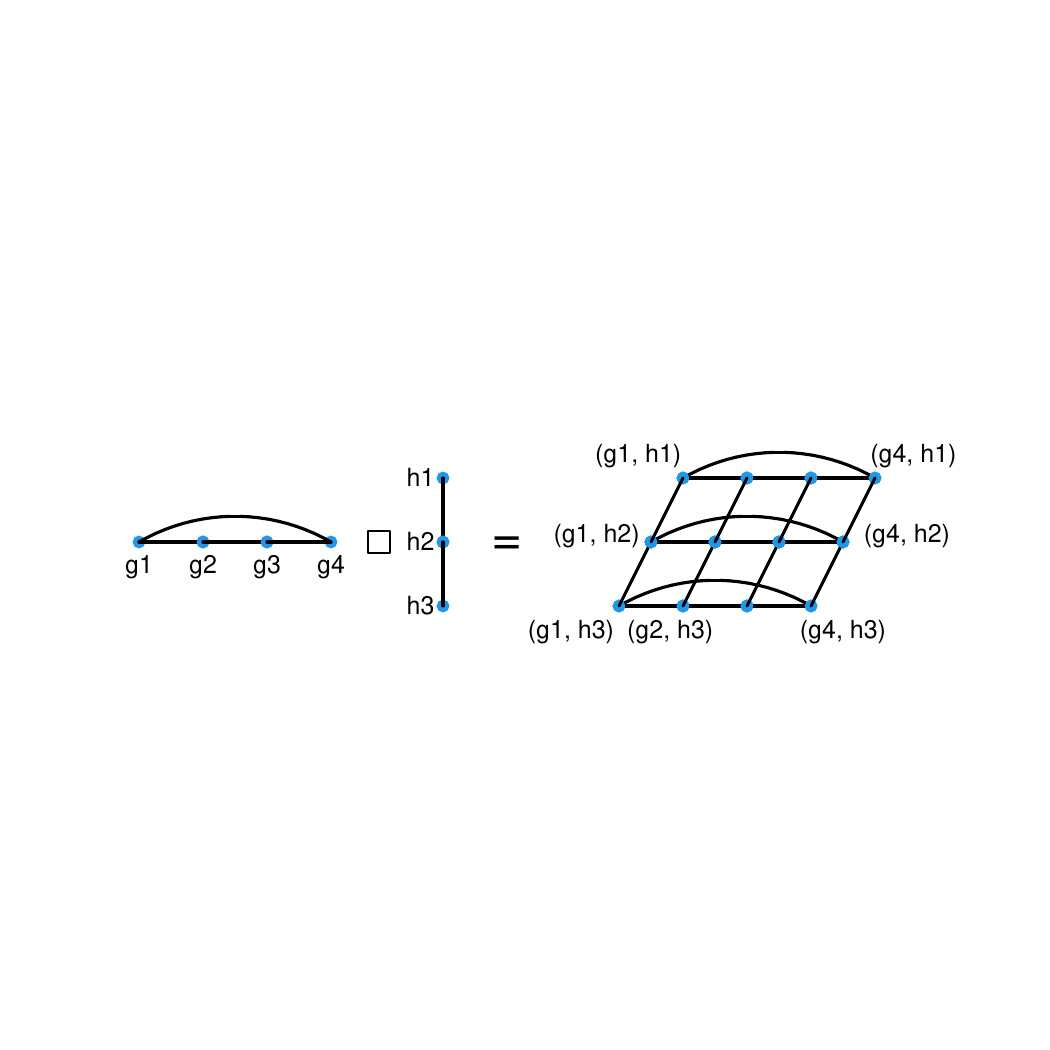} 
\end{tabular}
\vspace{-1.8in}
\end{center}
\caption{ Cartesian product graph $C_4 \Box P_3$, where $C_4$ is a
  cycle graph with 4 vertices and $P_3$ is a simple path graph with 3 vertices and 2 edges.}
\label{fig::cartgraph}
\end{figure}

\begin{figure}
  \begin{center}
\begin{tabular}{c}
    \includegraphics[width=4.5in]{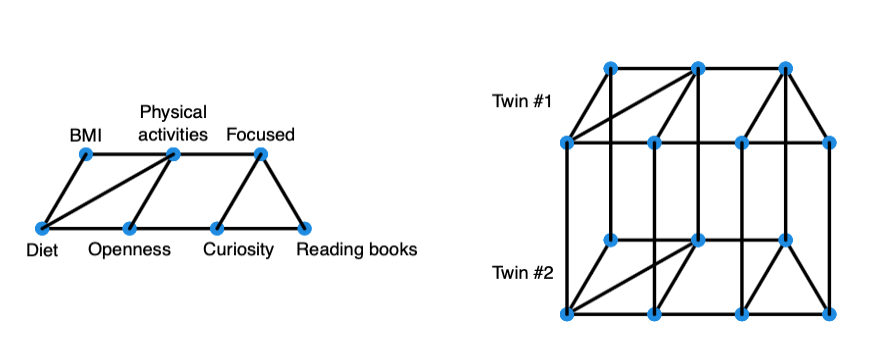}
\end{tabular}
 \vspace{-.2in}
  \end{center}
\caption{Cartesian product graph $G \Box K_2$, where $K_2$ is a complete graph with 2 vertices and 1 edge. Left panel: illustrative graph $G$ encodes
  the hypothetical conditional dependence relations among traits and hobbies as $V(G)$.
  Right panel: Prisms over graph $G$, formed by joining any vertex of $G$
  with its isomorphic image in $G'$; Only the same features are
  connected between the twins. }
\label{fig::carttwins}
\end{figure}

Along similar lines, the Bigraphical Lasso was proposed to
parsimoniously model conditional dependence relationships of matrix variate data based on the Cartesian
product of graphs~\citep{KLLZ13}.
The Cartesian product $G \Box H$ of graphs $G$ and $H$ (cf. Figure~\ref{fig::cartgraph})
is a graph such that  the vertex set is the Cartesian product $V(G) \times V(H)$ and 
two vertices $(g_1, h_1)$ and $(g_2, h_2)$ are adjacent in $G \Box H$
if and only if either $g_1 =g_2$ and $h_1$ is adjacent
to $h_2$ in $H$, or $h_1 = h_2$ and $g_1$ is adjacent to $g_2$ in $G$.
See Figure \ref{fig::carttwins} for illustration of the Cartesian 
product of graphs in modeling personality and behavior traits among twins.
A compelling justification for the proposed Kronecker sum model for
 the precision matrix is that similar models have been successfully
 used in fields including regularization of multivariate
 splines and design of physical
 networks; see~\cite{wood2006low} and \cite{imrich2008topics}.

\begin{figure}
 \begin{center}
\vspace*{-0.6in}
\begin{tabular}{l}
  \hspace*{-0.2in}
\includegraphics[width=0.98\textwidth,angle=0]{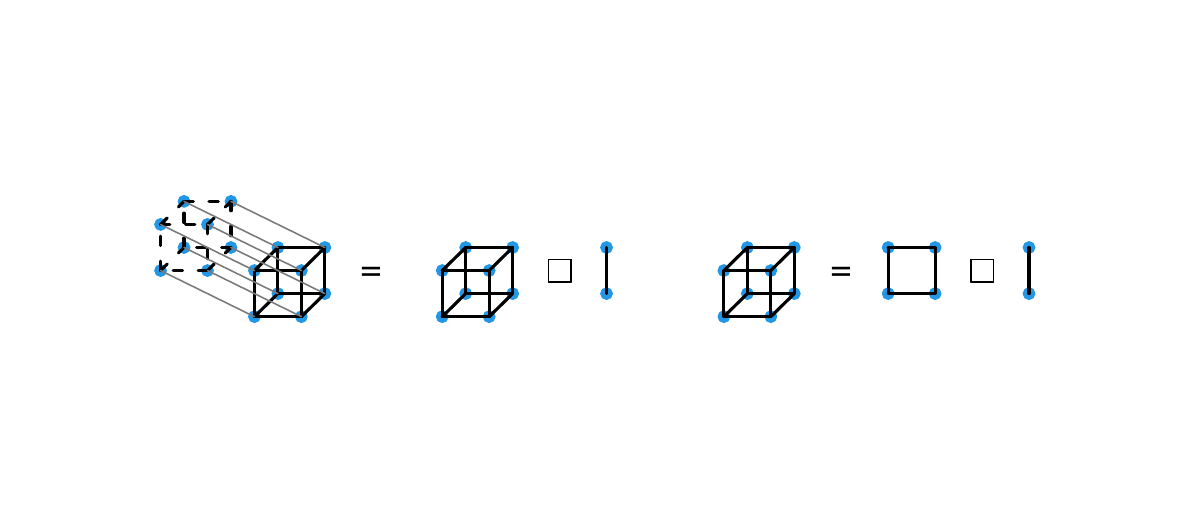}
\end{tabular}
\end{center}
\vspace*{-1.0in}
\caption{
  Cartesian product graph $K_2 \Box K_2 \Box K_2 \Box K_2$. The $n$-cube
  $Q_n, n \ge 1$ is defined as the $n^{th}$ power of $K_2$ with respect to the Cartesian product.}
\label{fig::cartgraph2}
\end{figure}

As pointed out by \cite{KLLZ13},  the associativity 
of the Kronecker sum naturally yields an approach to the modeling of datasets 
organized into 3 or higher-order tensors; cf. Figure~\ref{fig::cartgraph2}.
We demonstrate in~\cite{GZH19} that this model indeed generalizes
existing random matrix approaches to multilinear settings with more
than two axes of dependency structures well, by (a)
introducing a multiway tensor generalization of the Bigraphical 
Lasso estimator, known as the tensor graphical Lasso estimator,
for estimating sparse $L$-way decomposable inverse covariance matrices for
all integers $L \ge 2$; and (b) showing the rates of convergence in
the operator and Frobenius norm for estimating this class of inverse
covariance matrices for sub-gaussian tensor-valued data.
As a result, the  {\bf Te}nsor g{\bf ra}phical {\bf Lasso} (TeraLasso)
estimator is proven to effectively recover the conditional 
(in)dependence graphs and precision matrices for a class of Gaussian 
graphical models by restricting the topology to Cartesian product 
graphs; cf. Section~\ref{sec::tensordata}.


Consider the $L$-order random tensor $\X \in \mathbb{R}^{d_1
  \times \dots \times d_L}$, and assume that we are given $n$
independent samples $\X_1,  \ldots, \X_n \sim \X$. Here $\sim$
represents that two vectors follow the same distribution. Denote by $\mathbf{p} = [d_1,\dots,d_L]$ the vector of component dimensions and
$p$ the product of $d_j$s.
Hence
\ben
\label{eq::mk}
& \mathrm{vec}(\X) \in \mathbb{R}^p,  \quad \text{where } \quad
 p = \prod\nolimits_{k} d_k \quad \text{and} \quad m_k = \prod\nolimits_{i \neq k} d_i = {p}/{d_k}
\een
is the effective sample size we have to estimate the relations among the $d_k$ features
along the $k^{th}$ mode  in the tensor model.
It was shown in~\cite{GZH19} that due to the element replication inherent in the
Cartesian product structure, the precision matrix in the TeraLasso
model can be accurately estimated from limited data
samples of high dimensional variables with multiway coordinates such
as space, time and replicates.
Previously, we provided theoretical guarantees for the TeraLasso estimator~\eqref{eq::objfunc}, 
when the sample size is low, including single-sample convergence when 
$L \geq 3$~\citep{GZH19}. In particular, although single sample convergence was proved for $L >
2$,  empirically it was observed for all $L$.
In contrast, direct application of the models in \cite{FHT07} and the
analysis frameworks in \cite{RBLZ08}, \cite{ZLW08} and \cite{ZRXB11} require the sample size $n$ to scale proportionally to $p$,
which is still often too large to be practical. As a result, it is common to assume certain axes of $\X$ are i.i.d., often an overly simplistic model.

\subsection{ Contributions}
In the present work,  we strengthen the bounds for the relative errors in the operator and
Frobenius norm in~\cite{GZH19} by a factor of $\log p$,
improving upon those in Theorem~\ref{thm::orig}, as originally proved in~\cite{GZH19}.
These faster  rates of convergence are stated in
Theorem~\ref{thm::main} in the present paper.
We now show that the TeraLasso estimator achieves low errors with a 
{\em constant number} of replicates, namely $n = O(1)$, even for the 
$L=2$ regime.
This substantial improvement is due to the tighter error bounds on the
diagonal component of the loss function,
cf. Lemma~\ref{lemma::diagnew}.
This closes the gap between the finite (single) sample errors for the two-way models empirically 
observed in~\cite{GZH19} and the theoretical bounds therein.
The key technical innovation in the present work is the uniform concentration of measure bounds on the
trace terms appearing in the diagonal component of the loss
function~\eqref{eq::objfunc}, where we highlight tensor unfolding
techniques and Hanson-Wright inequalities.
Although the main results were presented in part in a conference
paper~\citep{ZG24}, we significantly expand the introduction to
illuminate the Kronecker Sum precision model, as well as provide the proof strategy and full proofs for the main theorems in Sections \ref{sec::proofofdiagfinal}, \ref{sec::strategy} and \ref{sec::proofofmainsupp}.

\subsection{Definitions and notations}
\label{sec::tensordata}
Let $e_1, \ldots, e_n$ be the canonical basis of $\R^n$.
Let $\Ball_2^n$ and $\Sp^{n-1}$ be the unit Euclidean ball and the 
unit sphere of $\R^n$, respectively.
For a set $J \subset \{1, \ldots, n\}$, denote
$E_J = \Span\{e_j: j \in J\}$.
We denote by $[n]$ the set $\{1, \ldots, n\}$.
We use $A$ for matrices, $\bA$ for tensors, and $\mathbf{a}$ for vectors.
For $\bA \in \mathbb{R}^{d_1 \times d_2 \ldots \times d_N}$, we use $\mathrm{vec}(\bA) \in \mathbb{R}^{d_1 \times d_2 \times \ldots   \times d_N}$ as in \cite{kolda2009tensor}, and
  define $\bA^T \in \mathbb{R}^{d_N\times \dots d_2 \times d_1}$ by analogy
  to the matrix transpose, i.e. $[\bA^T]_{i_1,\dots,i_N} = \bA_{i_N,\dots,   i_1}$.
The inner product of two tensors $\X, \Y \in \R^{d_1 \times 
  d_2 \times \ldots \times d_N}$ is sum of the products of their entries, i.e., 
\ben 
\label{eq::product}
\ip{\X, \Y} & = & 
\sum_{i_1=1}^{d_1} \sum_{i_{2}=1}^{d_{2}} \ldots\sum_{i_{N}=1}^{d_{N}}
x_{i_1 i_{2} \ldots \ldots i_{N}}  y_{i_1 i_2\ldots i_{N}}, 
\een 
where $x_{i_1,\dots,i_N}$ denotes the $(i_1,\dots, i_N)$-th element of
$\X$.
When extracted from the tensor, fibers are always assumed to be oriented as column vectors. 
The specific permutation of columns is not important so long as it is consistent across related calculations~\citep{kolda2009tensor}. 
Tensor unfolding of $\X$ along the $k$th mode is denoted as $\bX^{(k)}$, and is formed by
arranging the mode-$k$ fibers as columns of the resulting matrix of
dimension $d_k \times m_k$~\citep{kolda2009tensor}.
Denote by $\bX^{(k)T}$ its transpose. Denote by $X_j^{(k)}$ the $j^{th}$ column vector of $\bX^{(k)} \in
\R^{d_k \times m_k}$, where $d_k m_k  = p, \forall k \in [L]$.
One can compute the mode-$k$ Gram matrix $S^k$:
\ben
\label{eq::gramS2}
S^k = \bX^{(k)} \bX^{(k)T}/m_k
= \inv{m_k} \sum_{j=1}^{m_k} X_j^{(k)} \otimes X_j^{(k)} \in \R^{d_k \times d_k}.
\een
\subsection{The model and the method}
\label{sec::method}
For a subgaussian random variable $Z$, the $\psi_2$ norm of $Z$,
is defined as  $\norm{Z}_{\psi_2} = \inf\{t > 0\; : \; \E \exp(Z^2/t^2) \le 2 \}.$
\begin{definition}
  \label{def::vecZ}
Consider the tensor-valued data $\X$ generated from a subgaussian
random vector $Z = (Z_{j}) \in \R^p$ with independent mean-zero unit variance components whose $\psi_2$ norms are uniformly bounded:
\ben
\label{eq::tensordata}
\mvec{\X} & = & \Sigma_0^{1/2} Z, \; \; \text{ where } \; 
\expct{Z_{j}} = 0, \quad \E Z_{j}^2 =1,
\; \text{ and} \; \norm{Z_{j}}_{\psi_2} \leq C_0, \forall j.
\een
\end{definition}

We refer to $\X \in \mathbb{R}^{d_1 \times  \dots \times d_L}$ as an order-$L$ 
subgaussian random tensor with covariance $\Sigma_0 \in \mathbb{R}^{p
  \times p}$ for $\X$ as in~\eqref{eq::tensordata}.
Let $\X_1, \ldots, \X_n \in \mathbb{R}^{d_1 \times \dots \times d_L} \sim \X$ be $n$
i.i.d. random tensors following~\eqref{eq::tensordata}.
Let $\Omega_0 =\Sigma_0^{-1}$.
We assume that the precision matrix $\Omega_0 = \Psi_1 \oplus \dots \oplus \Psi_L $ of $\X$ is the $L$-way Kronecker sum of matrix components $\{\Psi_k\}_{k=1}^L$.
As such, we have
\ben
\label{eq::model}
\Omega_0  & = & \sum_{k =1}^L I_{[d_{1:k-1}]} \otimes \Psi_k \otimes I_{[d_{k+1:L}]}, 
\quad \text{where}  
\quad
{I}_{[d_{k:\ell}]} := \underbrace{{I}_{d_{k}} \otimes \dots \otimes {I}_{d_{\ell}}}_{\ell-k+1 \; \mathrm{factors}},
\een
where $\otimes$ denotes the Kronecker (direct) product and $\ell \geq
k$.
Denote by $S_n^{k}$ the mode-$k$ Gram matrix.
Now, we have $n m_k$ columns to compute the Gram matrices $S_n^{k},
\forall k$. 
Denote by $\Sigma_0^{(k)}, k \in [L]$ the corresponding factor-wise marginal 
covariance: $\Sigma_0^{(k)} = \mathbb{E}[S^{k}]$, for  $S^k$ as in
\eqref{eq::gramS2}. Then by linearity of expectations, 
\ben
\label{eq::gramK2}
&& S_n^{k} =\frac{1}{n m_k} \sum_{i=1}^n \bX^{(k,i)} [\bX^{(k,i)}]^T 
\quad \mathrm{and} 
\label{eq::gramE2}
\quad
\Sigma_0^{(k)} :=  \mathbb{E}[S_n^{k}]= \frac{1}{m_k}  \mathbb{E}[\bX^{(k)} \bX^{(k)T}].
\een
See~\cite{GZH19supp}.
The precision matrix  \eqref{eq::model}  has an immediate connection to the $L$ positive-semidefinite Gram
matrices $S_n^k \succeq 0  \in R^{d_k  \times d_k}$ associated with
each mode of the tensor $\X$, through tensor unfolding.
Denote by $\abs{\Omega}$ the determinant of $\Omega$.
Denote by $\mathcal{K}_{\mathbf{p}}^\sharp$ the set of positive
definite matrices that are decomposable into a Kronecker sum of
fixed factor dimensions $\mathbf{p} = [d_1, \ldots, d_L]$: 
\ben
\label{eq::kppintro}
\mathcal{K}_{\mathbf{p}}^\sharp
 & = & \{A \succ 0 | A \in \mathcal{K}_{\mathbf{p} } \subset
 \mathbb{R}^{p\times p}\}, \;\\
 \nonumber
 \text{ where}\;
 \; \mathcal{K}_{\mathbf{p}} & =& \{ {A}: \exists \: {B}_k \in \mathbb{R}^{d_k \times d_k} \: \mathrm{s.t.} \: {A} = {B}_1 \oplus \dots \oplus {B}_L\}.
\een
The TeraLasso estimator~\citep{GZH19} minimizes the negative
$\ell_1$-penalized Gaussian loglikelihood function $Q(\Omega)$
over the domain $\mathcal{K}_{\mathbf{p}}^{\sharp}$ of precision
matrices $\Omega \succ 0$,
where
\ben
\label{eq::lossfunc1}
Q(\Omega) & := & -\log \abs{\Omega} + \ip{\hat{S}, \Omega} +
\sum_{k=1}^L  m_k \rho_{n,k} \abs{{\Psi}_k}_{1,\off}, \; \text{ where
}  \\
\label{eq::gram}
\hat{S} & = & \inv{n} \sum_{i= 1}^n \mvec{\X_i^T} (\mvec{\X_i^T})^T,
\\
\nonumber
&& \; \text{ and} \; \forall k, \; \abs{\Psi_k}_{1,\off} = \sum_{i\neq j}{\abs{\Psi_{k,ij}}},
\een
and $\rho_{n,k}>0$ is a penalty parameter to be specified.
Here,  the objective function~\eqref{eq::lossfunc1} depends on the training 
data via the coordinate-wise Gram matrices $S_n^{k}$
\eqref{eq::gramK2} through projection, in view of \eqref{eq::model}, and 
the weight $m_k =p/d_k$ for each $k$ is determined by the number
of times for which a structure $\Psi_k$ is replicated in $\Omega_0$. 
This will become immediately obvious when we replace the trace term
$\ip{\hat{S},  \Omega}$ in~\eqref{eq::lossfunc1} with the weighted sum
over component-wise trace terms in~\eqref{eq::objfunc}; cf. Lemma~\ref{lemma::projection}. Then for $\mathcal{K}_{\mathbf{p}}^\sharp$ as in \eqref{eq::kppintro}, 
\ben
\label{eq::objfunc}
\lefteqn{\quad \quad    \quad \text{(TeraLasso)} \quad \hat{\Omega} 
:=\argmin_{\Omega \in \mathcal{K}_{\mathbf{p}}^\sharp} Q(\Omega) = } \\
\nonumber
&&  \!\!\!\!\!\!\!\!\!\!\!\!\!\argmin_{\Omega \in \mathcal{K}_{\mathbf{p}}^\sharp}
\big(-\log | \Omega| + \sum_{k=1}^L m_k \big(\langle S_n^{k},
\Psi_k \rangle + \rho_{n,k}  \abs{{\Psi}_k}_{1,\off} \big)\big).
\een
Here and in~\cite{GZH19}, the set of penalty parameters $\{\rho_{n,k}, k=1, \ldots, L\}$ are chosen to dominate the maximum of entrywise errors for estimating the population $\Sigma_{0}^{(k)}$~\eqref{eq::gramE2} with sample ${S}^{k}_{n}$ as in~\eqref{eq::gramK2}, for each $k \le L$ on event $\T$; cf. \eqref{eq::defineToffd}.
This choice works equally well for the subgaussian
model~\eqref{eq::tensordata}.

For $L=2$ and $\Omega_0 = \Psi_1\oplus \Psi_2=\Psi_1 \otimes
I_{d_2}+I_{d_1} \otimes \Psi_2$, the objective
function~\eqref{eq::objfunc} is similar in spirit to the BiGLasso
objective~\citep{KLLZ13}, where $S^k_n, k=1, 2$ correspond to the Gram
matrices computed from row and column vectors of matrix variate
samples $X_1, \ldots, X_n \in \R^{d_1 \times d_2}$ respectively.
When $\Omega_0 = \Psi_1 \otimes \Psi_2$ is a Kronecker product rather
than a Kronecker sum over the factors, the objective function
\eqref{eq::objfunc} is also closely related to the {\em Gemini
  estimators} by the first author of the present paper in~\cite{Zhou14a}, where $\log\abs{\Omega_0}$ is a linear combination of $\log\abs{\Psi_k}, k=1, 2$.
When $\X$ follows a multivariate Gaussian distribution and the precision matrix $\Omega_0$ has a decomposition of the form \eqref{eq::model}, the sparsity pattern of $\Psi_k$ for each $k$ corresponds to the conditional independence graph across the $k^{\text{th}}$ dimension of the data.
Similar to the graphical Lasso, incorporating an $\ell_1$-penalty promotes a sparse
graphical structure in the $\Psi_k$ and by extension $\hat{\Omega}$.
See for example \cite{ABG08,YL07,ZRXB11,RWRY08,KLLZ13, Zhou14a,GZH19,Horns19} and references therein.

\noindent{\bf More notation.}
We refer to a vector $x = (x_1, \ldots, x_n) \in \R^n$ with at most 
$d \in [n]$ nonzero entries as a $d$-sparse vector.
Denote by $\twonorm{x} = \sqrt{
  \sum_{i=1}^n x_i^2}$ and $\onenorm{x} := \sum_{j} \abs{x_j}$.
For a finite set $V$, the cardinality is denoted by $\abs{V}$.
For a given vector $x \in \R^n$, $\diag(x)$ denotes the diagonal matrix whose  
main diagonal entries are the entries of $x$.
For a symmetric matrix $A$, let $\phi_{\max}(A)$ and $\phi_{\min}(A)$ be the largest and the smallest  eigenvalue of $A$ respectively.
For a matrix $A$, we use $\twonorm{A}$ to denote its operator norm
and $\fnorm{A}$ the Frobenius norm, given by $\fnorm{A} = (\sum_{i, j} a_{ij}^2)^{1/2}$.
For a matrix $A = (a_{ij})$ of size $m \times n$,
let  $\norm{A}_{\infty} = \max_{i} \sum_{j=1}^n |a_{ij}|$ and  $\norm{A}_1 = \max_{j} \sum_{i=1}^m |a_{ij}|$ denote  the maximum absolute row and column sum of the matrix $A$ respectively.
Let $\norm{A}_{\max} = \max_{i,j} |a_{ij}|$.
Let $\diag(A)$ be the diagonal of $A$. Let $\offd(A) = A - \diag(A)$.
Let $\kappa(A) = \phi_{\max}(A)/\phi_{\min}(A)$ denote the condition
number for matrix $A$.
We use the inner product $\ip{A, B} = \tr(A^T B)$.
Fibers are the higher-order analogue of matrix rows and columns.
For two numbers $a, b$, $a \wedge b := \min(a, b)$, and 
$a \vee b := \max(a, b)$.
We write $a \asymp b$ if $ca \le b \le Ca$ for some positive absolute
constants $c,C$ that are independent of $n, m, p$,  and sparsity parameters.
Let $C, c, c', C_0, C_1, \ldots$ denote various absolute positive constants which may change line by line.

\noindent{\bf Organization.}
The rest of the paper is organized as follows.
Section~\ref{sec::theory} presents the main technical results, with discussions.
We elaborate on the new concentration of measure bounds
regarding the diagonal component of the loss function
in  Section~\ref{sec::proofofdiagfinal}, with full proof in 
Section~\ref{sec::appendeventMEproof}.
We conclude in Section~\ref{sec::conclude}.

\section{Theory}
\label{sec::theory}
In the present work, due to the tighter error bound on the diagonal
component of the loss function as stated in
Lemma~\ref{lemma::diagnew}, we achieve the sharper rates of
convergence in Theorem \ref{thm::main}, which significantly improve
upon earlier results in~\cite{GZH19} as stated in Theorem~\ref{thm::orig}.
Specifically, we replace the $p \log p$ in the earlier 
factor with $p$ for the relative errors in the operator and Frobenius norm in Theorem
\ref{thm::main} in the present work.
Under assumptions on the sparsity parameters, cf. Definition~\ref{def::A1} and
dimensions $d_k, \forall k \in [L]$, consistency and the rate of
convergence in the operator norm can be obtained for all $n$ and $L$.

\subsection{The projection perspective}
Throughout this paper, the subscript $n$ is omitted from $S_n^{k}$
and $\rho_{n,k}$ ($\delta_{n,k}$) in case $n=1$ to avoid clutter in the notation.
Lemma~\ref{lemma::projection} explains the smoothing ideas.
Intuitively, we use the $m_k$ fibers to estimate relations between and among the $d_k$ features along the $k^{th}$ mode, as encoded in $\Psi_k$.
Hence, this forms the aggregation of all data from modes other than $k$, which allows uniform concentration of measure bounds as shown in Lemma~\ref{lemma::diagnew} to be accomplished.

\begin{lemma}\textnormal{\bf (KS trace: Projection lemma)}
\label{lemma::projection}
Consider the mean zero $L$-order random tensor $\X \in \mathbb{R}^{d_1 \times \dots \times d_L}$.  Denote by $X_j^{(k)} \in \R^{d_k}$ the $j^{\text{th}}$ column vector in the matrix 
$\bX^{(k)} \in \R^{d_k \times m_k}$ formed by tensor unfolding.
Denote by $T := \ip{ \hat{S}, \Omega_0}$.
Then for sample covariance $\hat{S} := \mvec{\X^T} \otimes
\mvec{\X^T}$ and $\Omega_0$ as in~\eqref{eq::model}
\bens
\label{eq::D1}
T =
\sum_{k=1}^L \ip{m_k S^k, \Psi_k}
= \sum_{k=1}^L \sum_{j=1}^{m_k} \ip{\Psi_k, X_j^{(k)} \otimes X_j^{(k)} }, 
\eens
where $m_k S^k$ is the same as in \eqref{eq::gramS2}.
\silent{Then we also have
\bens
\label{eq::D2}
T =  \sum_{k=1}^L \sum_{i, j}^{d_k} \Psi_{k, ij} \ip{Y^{(k)}_i,   Y^{(k)}_j} 
 =  \sum_{k=1}^L \tr(\bY^{(k)} \Psi_k  \bY^{(k)T}) 
\eens}
Here $\mvec{A}$ of a matrix $A^{d_k \times m_k}$  is obtained by 
stacking columns of $A$ into a long vector of size $p = d_k 
\times m_k$.
\end{lemma}

\begin{lemma}
  \label{lemma::diagnew}
  Let $d_{\max} = \max_{k} d_k$ and $m_{\min}:=\min_{k} m_k$.
  Let $\Delta_{\Omega} \in   \mathcal{K}_{\mathbf{p}}$.
  Under the conditions in Lemma~\ref{lemma::projection},
  we have
  \bens
  \label{eq::diagfinal}
  \frac{ \abs{\ip{\diag(\Delta_{\Omega}), \hat{S} - \Sigma_0}}}{\twonorm{\Sigma_0}
    \fnorm{\diag(\Delta_{\Omega})}}
&\le & C_{\diag} \sqrt{d_{\max} L}
\big(1 + \sqrt{\frac{d_{\max}}{m_{\min}}}\big)
\eens
with probability at least $1-\sum_{k=1}^L 2 \exp(-c d_k)$.
\end{lemma}

\noindent{\bf Discussions.}
For simplicity, we state Lemma~\ref{lemma::projection} for the trace term $\ip{\hat{S}, \Omega_0}$ in case $n=1$, with obvious extensions for $n>1$ and for any $\Omega \in   \mathcal{K}_{\mathbf{p}}^\sharp$.
Now let $\bY^{(k)} = \bX^{(k)T}$. 
Denote by $Y_i^{(k)} \in \R^{m_k}$ the $i^{\text{th}}$ row vector in $\bX^{(k)}$. 
Then by~\eqref{eq::gramS2},
$$ Y_j^{(k)} \in \R^{m_k},  \forall j \in [d_k] \quad \text{ and } \quad
\forall i, j \in [d_k], m_k S^{k}_{ij}   = \ip{Y_{i}^{(k)}, Y_j^{(k)}},$$
which in turn can be interpreted as the tensor inner product~\eqref{eq::product} with $N = L-1$.
We mention in passing that $m_{\min}$ (resp. $n m_{\min}$) appears in 
the rates of  convergence in Theorems~\ref{thm::main} and~\ref{thm::orig}
as the effective sample size for estimating $\Omega_0$ for $n=1$
(resp. $n >1$). We discuss Lemma \ref{lemma::diagnew} further in
Section~\ref{sec::proofofdiagfinal}.
Before leaving this section, we define the support set of $\Omega_0$.
\begin{definition}\textnormal{\bf{(The support set of $\Omega_0$)}}
  \label{def::A1}
For each $\Psi_k$, $k = 1,\dots, L$, denote by 
$\supp(\offd(\Psi_k)) = \{(i,j): i \neq j, \Psi_{k, ij} \neq 0\}$. 
Let $s_k := \abs{\supp(\offd(\Psi_k))}$, for all $k$. 
Similarly, denote the support set of $\Omega_0$ by $\mathcal{S} = \{(i,j): i \neq 
j, \Omega_{0,ij} \neq 0 \}$, with $s:= \abs{\mathcal{S}}  = \sum_{k=1}^L m_k s_k$. 
\end{definition}

\subsection{The main results}
First, we state assumptions (A1), (A2) and (A3).
\bit
\item[(A1)]
  Let $\min_{k} m_k \ge \log p$.
  Denote $\delta_{n,k}  \asymp \twonorm{\Sigma_0} \sqrt{\frac{\log  p}{n 
    m_k}}$ for $k = 1, \ldots, L$.
  Let  $\rho_{n,k} = \delta_{n,k} /\ve_k$, where $0< \ve_k < 1$
   $\forall k$.

\item[(A2)]
  The smallest eigenvalue $\phi_{\min}({\Omega}_0)  = \sum_{k=1}^L \phi_{\min}({\Psi}_k) \geq\underline{k}_\Omega > 0$,
and the largest eigenvalue $\phi_{\max}({\Omega}_0) = \sum_{k=1}^L \phi_{\max}({\Psi}_k)  \leq\overline{k}_\Omega < \infty$.

\item[(A3)]
The sample size $n$ satisfies the following: for  some absolute
constant $C$,
\bens
\label{eq::mlower}
n (m_{\min})^2 \ge C^2 (L+1) \kappa(\Sigma_0)^4 (s \log p + L p),
\eens
where $m_{\min} := \min_{k} m_k$, $s = \sum_{k} m_k s_k $ is
as in Definition~\ref{def::A1}.
\eit

\begin{theorem}{\textnormal{\bf (Main result)}}
\label{thm::main}
  Suppose (A1), (A2), and (A3) hold. Then for absolute constants $C,
  c$, and  $C_L := C \sqrt{L+1}$, with probability $\geq 1  - L \exp(-c \log p)$, 
\bens
\nonumber
{\fnorm{\hat{\Omega} - \Omega_0}}/{\twonorm{\Omega_0}}
& \leq&  C \kappa(\Sigma_0) \big(\frac{s \log p + L p}{n m_{\min}}\big)^{1/2}, \\
  \label{eq::mainop}
{\twonorm{\hat{\Omega} - \Omega_0} }/{\twonorm{\Omega_0} }& \leq &
C_L \kappa(\Sigma_0) \big(\frac{{s \log p + L
      p}}{n m^2_{\min}}\big)^{1/2},  \\
\label{eq::mainfrob}
{\fnorm{\hat{\Omega} - \Omega_0}}/{\fnorm{\Omega_0}}
  & \leq &
 C_L \kappa(\Sigma_0) \big(\frac{{ s \log p + L p}}{n m^2_{\min}}\big)^{1/2}.
 \eens
\end{theorem}
The condition number for $\Sigma_0 = \Omega_0^{-1}$ is defined as
\bens
\kappa(\Sigma_0) =\kappa(\Omega_0) =
\twonorm{\Omega_0} \twonorm{\Omega_0^{-1}} =
 \frac{\sum_{k=1}^L 
   \phi_{\max}({\Psi}_k) }{\sum_{k=1}^L \phi_{\min}({\Psi}_k) },
 \eens
 where we have used the additivity of the eigenvalues of the Kronecker
 sum. Here and in~\cite{GZH19},
we focus on error bounds on the estimate of $\Omega_0$ itself, rather
than the individually factors.
We emphasize that we retain essentially the same error bound as that
in~\cite{GZH19}  for the off-diagonal component of the trace terms
in~\eqref{eq::objfunc}.
Event $\T$ is needed to control the off-diagonal component of the loss function:
\ben
\nonumber
\T & = & \bigcap_{k=1}^L \T_k \; \text{where }  \T_k =  \left\{ \max_{i \not=j} \abs{{S}^{k}_{n, ij} -
    \Sigma_{0,ij}^{(k)}} \le \delta_{n,k}\right\}, \\
\label{eq::defineToffd}
&& \text{ for } \delta_{n,k}  \asymp \twonorm{\Sigma_0}
\sqrt{{\log  p}/{(n m_k)}} > 0.
\een 
Intuitively, we use $nm_k$ fibers to estimate relations between and 
among the $d_k$ features along the $k^{th}$ mode as encoded in $\Psi_k$ and this allows optimal statistical rates of convergence to be derived, in terms of entrywise 
errors for estimating $\Sigma_{0}^{(k)}$ with ${S}^{k}_{n}$~\eqref{eq::gramK2}.
Correspondingly, events $\{\T_k, k=1, \ldots, L\}$ in
\eqref{eq::defineToffd}, which were originally defined in~\cite{GZH19}, cf. Proof of Lemma 12, are also used in the present work that
reflect this sample aggregation with $n m_k$ being the effective size
for estimating $\Psi_k$.

Indeed, as we will show in Theorem \ref{thm::orig}~\citep{GZH19},
these entrywise error bounds already enabled a significant improvement
in the sample size lower bound in order to estimate parameters and the
associated conditional independence graphs along coordinates
such as space, time and experimental conditions.
However,  these entrywise error bounds are not sufficient to achieve
the type of bounds
as in Theorem~\ref{thm::main} for inverse covariance estimation.
Using the entrywise error bounds to control the diagonal components of
the trace terms will result in an extra $\log p$ factor in the sample size
lower bound and correspondingly a slower rate of convergence. This
extraneous $\log p$ factor is undesirable since the diagonal component
of the loss function dominates the overall rate of convergence in
sparse settings for inverse covariance estimation.

\noindent{\bf Summary.}
    The worst aspect ratio is defined as 
    $$\max_{k} ({d_k}/{m_k}) = \frac{d_{\max}}{m_{\min}}=\frac{p}{m_{\min}^2}.$$
    Clearly, a smaller aspect ratio implies a faster rate of convergence 
    for the relative errors in the operator and Frobenius norm.
First, observe that for relative error in the operator norm in  Theorem~\ref{thm::orig},
$(L+1)(s+p) \log p$ therein is replaced with $s \log p + Lp$ cf.
Theorem~\ref{thm::main}.
The same improvement holds true for the Frobenius norm error.
Here we eliminate the extraneous $\log p$ factor from the diagonal
component of the error through new concentration of measure analysis
in the present work; cf. Lemma~\ref{lemma::diagnew}.
This is a significant improvement for two reasons: (a) since $p$
is the product of the $d_k$s, $\log p = O(\sum_{k} \log d_k)$ is often
nontrivial, especially for larger $L$; and (b) more importantly, 
for $L=2$ and $n = O(1)$  (in contrast to $L > 2$), the error bound in
the operator norm in Theorem \ref{thm::orig} by~\cite{GZH19} will \emph{diverge} for any $s \geq 0$ as $p = d_1 d_2$ increases, since
\ben
\label{eq::K2imp}
&& \frac{p \log p}{m_{\min}^2} = \frac{d_1 d_2 \log p}{(d_1 \wedge d_2)^2} \ge
\log p,
\een
where $m_{\min} = p/(d_1 \vee d_2) = d_1 \wedge d_2$
and equality holds only when $d_1=d_2$. 
As a result, in Theorem \ref{thm::orig}~\citep{GZH19}, the sample lower
bound, namely, $n(m_{\min})^2 \ge C^2 \kappa(\Sigma_0)^4 (s+ p)
(L+1)^2 \log p$ implies that $n = \Omega(\log p)$, since $m_{\min}^2
\le p$ in view of~\eqref{eq::K2imp}. In contrast, the lower bound
on $n m^2_{\min}$  in (A3) is less stringent, saving a factor of
$O(\log p)$.

This is consistent with the successful finite sample experiments in~\cite{GZH19},
where for $L=2$, bounded errors in the operator norm are observed as $p$ increases.
As a result, our new bound supports the use of the TeraLasso estimator
when $L=2$, so long as {\em a small number of replicates} are
available, that is, when $n = o(\log p)$, in a way that the previous
Theorem~\ref{thm::orig} cannot. More precisely, for {\em finite sample
  settings}, namely, when $n = O(1)$, the relative errors will
still be bounded at $O_p(1)$ for $L=2$, for example, when the two dimensions
are at the same order: $d_1  \asymp d_2$, and rapidly converge to zero for $L >
2$; cf.~Theorem~\ref{thm::maincubic}.

\noindent{\bf Single sample convergence.}
First of all, both Theorems~\ref{thm::main} and~\ref{thm::orig} imply $n=1$
convergence for the relative error in the operator norm, when $L \ge 3$
and $d_1 \asymp \ldots \asymp d_L$, which we refer to as the cubic tensor
settings, since potentially $m^2_{\min} \ge m_{\min}  d_{\max} \log p
= p \log p$ will hold.
However, when the $d_k$s are skewed, this may
not be the case.  To make this clear, we first state Corollary \ref{thm::cubic}.
\begin{corollary}[Dependence on aspect ratio for $n=1$]
  \label{thm::cubic}
    Suppose (A1), (A2) and (A3) hold for $n=1$.
  Then  with probability at least $1-L \exp(c \log p)$, we have for some absolute constants $c, C$,
  \bens
\frac{\twonorm{\hat{\Omega} - \Omega_0}}{\twonorm{\Omega_0}}
\vee \frac{\fnorm{\hat{\Omega} - \Omega_0}}{\fnorm{\Omega_0}}
\le  C \kappa(\Sigma_0) \big({\frac{L d_{\max}}{m_{\min}}}\big)^{1/2} \big(\sum_{k=1}^L \frac{s_k \log p}{d_k} + L\big)^{1/2}.
\eens
\end{corollary}
Under the bounded aspect ratio regime, 
the relative errors in the operator and Frobenius norm for estimating 
the precision matrix $\Omega_0$ depend on the decay of the worst  aspect 
ratio ${d_{\max}}/{m_{\min}}$ and the average of ${s_k \log p}/{d_k}$
over all modes, which represents relative sparsity levels (sparsity / dimension) in an 
average sense.
For $L > 2$, typically the aspect ratio is much less than 1 and
convergence happens rapidly. If the sparse support set is small
relative to nominal dimension $d_k$ along each mode,
for example, when $\frac{s_k \log p}{d_k} = O(1)$, this
convergence is at the rate of decay of the worst aspect ratio.
 In this case, the diagonal component dominates the rate of convergence and this is 
essentially optimal, since in the largest component with dimension
$d_{\max}$, it has $d_{\max}$ parameters to be estimated
and $m_{\min} = p/d_{\max}$ effective samples for the task.
Moreover, $L$ is needed in the bound since we estimate $L$
components all together  using one sample in case $n=1$.

\subsection{Cubic tensor and optimality}
\label{sec::cubic}
 As a final example, we consider the \textbf{cubic} setting,
 where $d_1 \asymp \ldots \asymp d_L \asymp p^{1/L}$.
In words, a tensor is cubical if all $d_j$s are at the same order.
Then
\ben 
\label{eq::dmratio}
\text{aspect ratio} \quad :=\frac{d_{\max}}{m_{\min}} \asymp \frac{p^{1/L}}{p^{1-1/L}} =p^{2/L -1}.
\een
Note that for $L > 2$, we obtain a fast rate of 
convergence in the operator norm for $n=1$, since in the cubic 
tensor settings, the effective sample size $m_{\min}$ increases
significantly faster than $\sqrt{p}$ given that $d_{\max} = o(p^{1/2})$.
More precisely, we state Theorem~\ref{thm::maincubic}, where we
consider the cubic tensor setting and $n=1$.
\begin{theorem}{\bf{(The cubic tensor)}}
  \label{thm::maincubic}
Under the conditions in Theorem~\ref{thm::main}, suppose $d_k =O(m_k)$ for all $k$.
Suppose $m_1 \asymp m_2 \asymp \ldots \asymp m_L$.  Then,
 \bens 
&& \!\!\!\!\!\!\!\!\!\!\frac{\fnorm{\hat{\Omega} - \Omega_0}}{\kappa(\Sigma_0)\twonorm{\Omega_0} }
= O_P\big(\big(\sum_{k=1}^L s_k \log p + L d_{\max}\big)^{1/2}
\big), \text{ and} \\
&&
\label{eq::densecubic}
 \!\!\!\!\!\!\!\!\!\!\frac{\twonorm{\hat{\Omega} -
     \Omega_0}}{\kappa(\Sigma_0)\twonorm{\Omega_0}}  = O_P\big(
 \big(L\sum_{k=1}^L s_k \log p + L^2 d_{\max}\big)^{1/2}/{m_{\min}^{1/2}}  \big).
\eens
Suppose in addition 
$d_1 \asymp \ldots \asymp d_L = \Omega\big(({\log p}/{L}) \sum_{k} s_k
\big).$ Then
\bens
\label{eq::cubicopt}
\frac{\twonorm{\hat{\Omega} - \Omega_0}}{\twonorm{\Omega_0}} \vee \frac{\fnorm{\hat{\Omega} - \Omega_0}}{\fnorm{\Omega_0}} = O_P \big(L \kappa(\Sigma_0) p^{{1}/{L} - {1}/{2}}\big).
\eens
\end{theorem}
Theorem~\ref{thm::maincubic}
shows that convergence will occur for the {\em dense cubic}
case, so long as $m_{\min} =p/d_{\max}= \Omega\big(L \log
p \sum_{j=1}^L s_j  + L^2  d_{\max}\big),$
which is a reasonable assumption in case $L > 2$ and holds under
(A3). In other words, the relative errors in the operator and
Frobenius norm are bounded so long as the effective sample size $m_{\min}$ is at least $L^2 d_{\max} \ge L \sum_{k} d_{k}$, which is roughly $L$ times the total number of (unique) diagonal entries in $\{\Psi_k, k =1, \ldots, L \}$, and also at least  $L
     \log p$ times $\sum_k s_k$, which in turn denotes the size of
     total supports $\sum_k  \abs{\mathcal{S}_k}$ over off-diagonal components of
     factor matrices $\{\Psi_1, \ldots,    \Psi_k\}$.
 Consider now an even more special case.
 Suppose that in the cubic tensor setting, we have
 $d_{\max} =\Omega( \log p  \sum_{j} s_j/L)$  in addition.
Then the error in the operator norm is again dominated by the square root of
 the aspect ratio parameter.
In other words, to achieve the near optimal rate of $O_P(p^{{1}/{L} - {1}/{2}})$,
it is sufficient for each axis dimension $d_k, k \in [L]$ to dominate
the \emph{average sparsity} across all  factors, namely, $\sum_{k}
s_k/L$ by a $\log p$ factor. A more general result has been stated in
Corollary~\ref{thm::cubic}. The proof of Theorem~\ref{thm::main}
appears in Section~\ref{sec::strategy}. We prove Theorem~\ref{thm::maincubic} and 
Corollary \ref{thm::cubic} in Sections~\ref{sec::proofmaincubic}
and~\ref{sec::proofofcubic} respectively.

\subsection{Related work}
Models similar to the Kronecker sum precision model have been
successfully used in a variety of fields, including regularization of
multivariate
splines~\citep{wood2006low,eilers2003multivariate,kotzagiannidis2017splines,wood2016},
design of physical networks
\citep{imrich2008topics,van2000ubiquitous,fey2018splinecnn},
neuroscience~\citep{GPZN17}, and Sylvester equations arising from the
discretization of separable
$L$-dimensional PDEs with tensorized finite
elements~\citep{grasedyck2004existence,kressner2010krylov,beckermann2013error,shi2013backward,ellner1986new}.
Additionally, Kronecker sums find extensive use in applied mathematics and
statistics, including beam propagation physics
\citep{andrianov1997matrix}, control
theory~\citep{luenberger1966observers,chapman2014controllability},
fluid dynamics~\citep{dorr1970direct},
errors-in-variables~\citep{RZ17}, and spatio-temporal modeling and
neural processes \citep{schmitt2001numerical,GPZ17,fan19a}.
When the data indeed follows a matrix normal model, the
BiGLasso~\citep{KLLZ13} and TeraLasso~\citep{GZH19} also effectively
recover the conditional dependence graphs and precision matrices
simultaneously for a class of Gaussian graphical models by restricting
the topology to Cartesian product graphs.
We provided a composite gradient-based optimization algorithm,
and obtained algorithmic and statistical rates of convergence for
estimating structured precision matrix for tensor-valued
data~\citep{GZH19}.

Recently, several methods have arisen that can speed up the numerical
convergence of the optimization of the BiGLasso objective of
\cite{KLLZ13}, cf.~\eqref{eq::objfunc} with $L=2$.
A Newton-based optimization algorithm for $L=2$ was presented in
\cite{yoon2021eiglasso} that provides significantly faster convergence
in ill-conditioned settings.
Subsequently, \cite{li2022two} developed a scalable flip-flop
approach, building upon the original BiGLasso flip-flop algorithm as
derived in \cite{KLLZ13}. Using the Kronecker sum eigenvalue decomposition similar to that of
\cite{GZH19} to make the memory requirements scalable, their algorithm
also provides faster numerical convergence than the first-order
algorithm presented in \cite{GZH19}.
They also provided a Gaussian copula approach for applying the model to certain non-Gaussian data.
Subsequent to~\cite{GZH19}, a related SG-PALM was presented in~\cite{wang2021sg}, where the precision matrix is the \emph{square} of an $L$-way Kronecker sum.
See \cite{wang2022kronecker} for a survey of multiway covariance models.

As mentioned, normality is not needed in our proofs; instead, we
consider subgaussian ensembles and derive tight concentration of
measure bounds, using tensor unfolding techniques. For recent
concentration of measure results on subgaussian matrix-variate models, we refer
to~\cite{RZ17},~\cite{Zhou19}, and~\cite{Zhou24}.

\section{The new concentration bounds}
\label{sec::proofofdiagfinal}
Throughout this proof, we assume $n=1$ for simplicity.
We now provide outline for proving the upper bound on the diagonal
component of the main result of the paper. Recall the true parameter $\Omega_0 =\Psi_1
\oplus \dots \oplus \Psi_L$, where $\Psi_k \in \R^{d_k \times
  d_k}$~\eqref{eq::model}.
Since $\Omega_0 \in \ksump$, we have
\ben
\label{eq::DO}
\forall \Omega \in \ksump, \; \;
\Delta_\Omega  := \Omega - \Omega_0 =  \Delta_{\Psi_1}
\oplus \Delta_{\Psi_2}  \oplus \ldots \oplus \Delta_{\Psi_L},
\een
for some $\Delta_{\Psi_k} \in \R^{d_k \times d_k}$ whose off-diagonal (but not diagonal)
elements are uniquely determined.
For self-containment, we state Lemma~\ref{lemma::orthoDelta},
where we also state the notation we use throughout
this section. Here we use the trace-zero convention which guarantees 
the uniqueness of the $\Delta'_{\Psi_k}$ in \eqref{eq::original}.
We will then restate Lemma~\ref{lemma::diagnew} in 
Lemma~\ref{lemma::diagfinal}.
The off-diagonal component has been dealt with in~\cite{GZH19supp};
cf. Lemmas 11 and 12 therein.
Proof of Lemmas is deferred to Section~\ref{sec::appendeventMEproof}.

\begin{lemma}\textnormal{{\bf (Decomposition lemma)}~\citep{GZH19supp}}
  \label{lemma::orthoDelta}
  Let $\Omega \in \ksump$.
  Then $\Delta_{\Omega} = \Omega - \Omega_0 \in \ksump$.
  To obtain a uniquely determined representation,
we rewrite~\eqref{eq::DO} as follows: 
\ben
\nonumber
\Delta_\Omega 
&  = &
\Delta'_\Omega +\tau_{\Omega} I_p, \quad \text{ where} \quad
\tau_{\Omega} = \tr(\Delta_\Omega)/p,  \quad \text{ and } \\
\label{eq::original}
&& \Delta'_\Omega = 
\Delta'_{\Psi_1}  \oplus \ldots \oplus \Delta'_{\Psi_L}, \; \text{
  where } \; \; \tr(\Delta'_{\Psi_k})= 0  \; \text{for all $k$}.
\een
Thus we have
 \ben
 \diag(\Delta'_\Omega) 
& = &
\label{eq::diagsum}
\label{eq::diagDK22}
\sum_{k=1}^L \diag( \tilde{\Delta}_k) \quad \text{ where} \; \;
\diag(\tilde{\Delta}_k)  :=  I_{[d_{1:k-1}]} \otimes \diag(\Delta'_{\Psi_k}) \otimes
I_{[d_{k+1:L}]},
\een
and moreover,
\ben 
\label{eq::ortho}
&& \quad \fnorm{\diag(\Delta_\Omega)}^2 =
\sum_{k=1}^L m_k \fnorm{\diag(\Delta'_{\Psi_k})}^2  + p
\tau_{\Omega}^2, \\
\nonumber
\label{eq::opnormbound}
&& 
\sum_{k=1}^L \sqrt{d_{k}}
\fnorm{\diag(\Delta'_{\Psi_k})} \le \sqrt{\frac{L d_{\max}}{m_{\min}}} \fnorm{\diag(\Delta_{\Omega})}.
\een
\end{lemma}

\begin{proof}
  The existence of such parameterization in \eqref{eq::original} is given in Lemma 7
  \cite{GZH19supp}, from which \eqref{eq::ortho} immediately follows, by  orthogonality of the decomposition.
  Now we have by elementary inequalities:
\bens 
&& \sum_{k=1}^L \sqrt{d_k} \fnorm{\diag(\Delta'_{\Psi_k})} =
\sum_{k=1}^L \sqrt{\frac{d_k}{m_k}}  \sqrt{m_k} \fnorm{\diag(\Delta'_{\Psi_k})} \\
&&\le
\max_{k} \sqrt{\frac{d_k}{m_k}}  \sqrt{L}\big(\sum_{k=1}^L m_k 
\fnorm{\diag(\Delta'_{\Psi_k})}^2\big)^{1/2}  \le 
\frac{\sqrt{d_{\max}}}{\sqrt{m_{\min}}}
\sqrt{L} \fnorm{\diag(\Delta_\Omega)}.
\eens 
Thus the lemma holds in view of~\eqref{eq::ortho}.
\end{proof}

\begin{lemma}{\textnormal{\bf (New diagonal bound)}}
  \label{lemma::diagfinal}
Following the notation as in Lemma~\ref{lemma::orthoDelta},   where $\tau_{\Omega} = \tr(\Delta_\Omega)/p$, we have
with probability at least $1 - \sum_{k} \exp(-c d_k) -c'/{p^4}$, 
\bens
\abs{\ip{\diag(\Delta_{\Omega}), \hat{S} - \Sigma_0}
}/\twonorm{\Sigma_0} \le  C_0 \sum_{k=1}^L d_k  \fnorm{\diag(\Delta'_{\Psi_k})} +  C_1 \sqrt{L d_{\max}} \fnorm{\diag(\Delta_{\Omega})},
\eens
where $c, c', C_0, C_1$ are absolute constants, and hence  Lemma~\ref{lemma::diagnew} holds.
\end{lemma}

Note when we have $d_k  \le \sqrt{p}$
for all $k$, or equivalently, when $\max_{k} \sqrt{\frac{d_k}{m_k}}
\le 1$, we do not need to pay the extra factor of $\sqrt{\log p}$ as in 
Lemma 13~\citep{GZH19supp} on the diagonal portion of the error bound,
resulting in the improved rates of convergence in Theorem~\ref{thm::main}. Note that when $d_k
=o(m_k \log p) , \forall k$, the bound in Lemma~\ref{lemma::diagfinal}
still leads to an improvement on the overall rate.

\begin{lemma}
\label{lemma::eventME}
Let $\Sp^{d_k-1}$ be the sphere in $\R^{d_k}$.
Construct an $\ve$-net $\Pi_{d_k} \subset \Sp^{d_k-1}$
such that $\abs{\Pi_{d_k}} \le (1 + 2/\ve)^{d_k}$, where $0< \ve <
1/2$, as in Lemma~\ref{lemma::net}. Recall $\bY^{(k)} = (\bX^{(k)})^T$.
Let $\delta  =(\delta_1, \ldots, \delta_{d_k})$.
Let $C_m, c$ be some   absolute constants.
Define the event $\MG_k$ as:
\ben
  \label{eq::tk}
&& \sup_{\delta \in \Pi_{d_k}}
 \sum_{i=1}^{d_k} \delta_{i} \big(\ip{Y^{(k)}_i,  Y^{(k)}_i}
- \E \ip{Y^{(k)}_i,  Y^{(k)}_i}\big) \le
t_{k}, \\
&& \text{where} 
\quad  t_{k}   := C_m   \twonorm{\Sigma_0}  (\sqrt{p} \vee d_k ).
\een
Let $\MG = \MG_1 \cap \ldots \cap \MG_L$.
  Then $\prob{\MG} \ge 1 -\sum_{k} \exp(-c d_k)$.
Moreover, we have by a standard approximation argument,
on event $\MG$,
\bens
\text{ simultaneously for all } k, \; \; 
\sup_{\delta \in \Sp^{d_k-1}} \sum_{i=1}^{d_k}
\delta_{i} \big(\ip{Y^{(k)}_i, Y^{(k)}_i} - \E \ip{Y^{(k)}_i,
    Y^{(k)}_i} \big) \le \frac{t_{k}}{1-\ve}.
\eens
\end{lemma}

\noindent{\bf Proof idea.}
Notice that the expression for $t_{k}$ clearly depends on the 
dimension $d_k$ of  $\Psi_k$. Let $\delta \in \mathbb{R}^{d_k}$.
Using the notation in Lemma~\ref{lemma::orthoDelta},
let $\diag(\Delta'_{\Psi_k}) = \diag(\delta_{1}, \ldots,
\delta_{d_k})$ and
  \ben
  \label{eq::diagDK}
  \diag(\tilde{\Delta}_k) :=
  I_{[d_{1:k-1}]} \otimes \diag(\Delta'_{\Psi_k}) \otimes 
  I_{[d_{k+1:L}]}.
  \een
  Now for each $1 \le k \le L$,
following Lemma~\ref{lemma::projection}, we have
\ben
\nonumber 
\lefteqn{\ip{\diag( \tilde{\Delta}_k), \hat{S} - \Sigma_0} =
m_k \ip{S^k - \E (S^k), \diag(\Delta'_{\Psi_k}) }}\\
\nonumber
&= &
\tr(\bY^{(k)} \diag(\Delta'_{\Psi_k})  \bY^{(k)T}) - \E \tr(\bY^{(k)}
\diag(\Delta'_{\Psi_k})  \bY^{(k)T}) \\
& &
\label{eq::YYtrace}
\quad \quad = \sum_{j=1}^{d_k} \delta_{j} \big(\ip{Y^{(k)}_j,  Y^{(k)}_j} - \E 
  \ip{Y^{(k)}_j,  Y^{(k)}_j}\big). 
\een
To bound the probability for event $\MG_k$, we use the Hanson-Wright
inequality in \cite{RV13}, cf. Theorem 1.1 therein, and the union bound.
The rest is deferred to Section~\ref{sec::eventME}.

\section{Proof of Lemmas~\ref{lemma::diagfinal}
and~\ref{lemma::eventME}}
\label{sec::appendeventMEproof}
Let the sample covariance $\hat{S} := \mvec{\X^T} \otimes \mvec{\X^T}$
be as in \eqref{eq::gram} and  $\Sigma_0 = \Omega_0^{-1}\in \mathbb{R}^{n
  \times n}$ be the true covariance matrix.
Let $Z \in \mathbb{R}^{p}$ denote an 
isotropic sub-gaussian random vector with independent coordinates as
in Definition~\ref{def::vecZ}.
Let
\ben
\label{eq::diagDK2}
\diag(\tilde{\Delta}_k) & := & I_{[d_{1:k-1}]} \otimes \diag(\Delta'_{\Psi_k}) \otimes 
I_{[d_{k+1:L}]}.
\een
This explains~\eqref{eq::WZ}. Consequently, by \eqref{eq::diagDK2}
$$\fnorm{\diag(\tilde{\Delta}_k)}^2 :=  m_k
\fnorm{\diag(\Delta'_{\Psi_k})}^2.$$
See also \eqref{eq::objfunc}.
Indeed, as expected, $\tr(\hat{S})$ converges to 
$\tr(\Sigma_0)$ at the rate of 
\bens 
\abs{\tr(\hat{S}) - \tr(\Sigma_0)}/p =O_P( \twonorm{\Sigma_0} 
\sqrt{\log p/(np)}). 
\eens

\silent{
To bound $\drp_0$,  we rewrite
\bens 
\ip{\hat{S} - \Sigma_0, I} = \tr(\hat{S} - \Sigma_0) = Z^T \Sigma_0 Z -\E (Z^T \Sigma_0 Z), 
\eens 
where $\Sigma_{0} \succ 0$ is a $p \times p$ symmetric positive definite 
matrix and $Z \in \R^p$ is the same as in~\eqref{eq::tensordata}.
Set $t = C \sqrt{p \log p} \twonorm{\Sigma_0 }$.
Thus, we have by the Hanson-Wright inequality,
\bens
\lefteqn{\prob{\abs{\ip{\hat{S} - \Sigma_0, I} > C  \twonorm{\Sigma_0 }
      \sqrt{p \log p}} }} \\
& = & 
\prob{\abs{Z^T \Sigma_0  Z - E{Z^T \Sigma_0  Z}} > C \sqrt{p \log p}
  \twonorm{\Sigma_0 }}\\ 
& \le & 
2 \exp\left(-c\min\left(\frac{C^2 p \log p \twonorm{\Sigma_0}^2}{\fnorm{\Sigma_0 }^2}, 
    C \sqrt{p \log p}\right) \right)  \le  \inv{p^4},
\eens
where
$(C^2 \wedge C) c \ge
4$ and $\fnorm{\Sigma_{0} } \le \sqrt{p} \twonorm{\Sigma_0 }$.
Hence by the union bound, we have by Lemma~\ref{lemma::eventME} and 
the bound on $\drp_0$ immediately above, 
\bens 
\prob{\MG \cap \drp_0} \ge 1- c \exp(-\log p) - \sum_{k} \exp(-c d_k).
\eens
The lemma thus holds upon adjusting the constants. }

First we show the following bounds on the $\ve$-net of $\Sp^{d_k-1},
\forall k$. 
\begin{lemma}{\textnormal{~\cite{MS86}}}
  \label{lemma::net}
Let $1/2> \ve > 0$. For each $k \in [L]$, one can construct an
$\ve$-net $\Pi_{d_k}$, which satisfies $$\Pi_{d_k} \subset \Sp^{d_k-1} \; \text{ and }\; \abs{\Pi_{d_k}}  \le (1+2/\ve)^{d_k}.$$
\end{lemma}

By Lemma~\ref{lemma::projection}, we have for the diagonal and
off-diagonal components of the trace term defined as follows: for
$\Omega_0 =   \Psi_1 \oplus \dots 
  \oplus \Psi_L$, 
\bens
\label{eq::L2.7}
\ip{ \hat{S},\diag(\Omega_0)} 
& = &
\sum_{k=1}^L  
\sum_{i=1}^{d_k} \Psi_{k, ii} \ip{Y^{(k)}_i,  Y^{(k)}_i} \quad \text{ and } \\
\nonumber
\ip{ \hat{S},  \offd(\Omega_0)}
& = &
\sum_{k=1}^L \sum_{i\not=j}^{d_k} \Psi_{k,
ij} \ip{Y^{(k)}_i,  Y^{(k)}_j},
\eens
where $\diag(\Omega_0) = \diag(\Psi_1) \oplus \dots 
\oplus \diag(\Psi_L)$ and $\offd(\Omega_0) = \offd(\Psi_1) \oplus \dots 
\oplus \offd(\Psi_L)$. See~\eqref{eq::YYtrace}, for which such a
decomposition is useful.

\subsection{Proof of Lemma~\ref{lemma::diagfinal}}
\begin{proofof2}
Besides $\MG$,
we need the following event $\drp_0$: 
\ben
\label{eq::diagsumproof}
\quad \drp_0 = \left\{\abs{\ip{I_p, \hat{S} - \Sigma_0} }
  \le   C \sqrt{p \log p} \twonorm{\Sigma_0}\right\}.
\een
Suppose $\MG \cap \drp_0$ holds. Denote by
\ben
  \label{eq::diagDK3}
\diag(\Delta'_{\Psi_k}) =  \diag(\delta^k_{1}, \ldots,
\delta^k_{d_k}) =: \diag(\delta^k),
\een
where $\twonorm{\delta^k} :=  \fnorm{\diag(\Delta'_{\Psi_k})}$.
Denote by
\bens
\label{eq::defineTKprime}
t_k' & = &  t_{k}  \twonorm{\delta^k}=
C_m  \twonorm{\Sigma_0} \fnorm{\diag(\Delta'_{\Psi_k})} \big(\sqrt{p} \vee d_k \big),
\eens
for $t_{k}$ as in \eqref{eq::tk}.
For each index $1 \le k \le L$, on event $\MG_k$, simultaneously
for all $\diag( \tilde{\Delta}_k)$ as in \eqref{eq::diagsum} and~\eqref{eq::diagDK2}, we have
\bens
\label{eq::Wquad}
\lefteqn{\abs{\ip{\diag(\tilde{\Delta}_k), \hat{S} - \Sigma_0}}
 = \abs{m_k \ip{S^k - \E S^k, \diag(\Delta'_{\Psi_k})} }}\\
 \nonumber 
& = &
 \nonumber 
\twonorm{\delta^k}  \abs{\sum_{j=1}^{d_k} \frac{\delta^k_{j}}{\twonorm{\delta^k}}
   \big(\ip{Y^{(k)}_j,  Y^{(k)}_j} - \E 
\ip{Y^{(k)}_j,  Y^{(k)}_j}\big)} \\
& \le &
\twonorm{\delta^k}
\sup_{\delta \in \Sp^{d_k-1}} \sum_{i=1}^{d_k} \delta_{i} \big(\ip{Y^{(k)}_i, Y^{(k)}_i} -  \E \ip{Y^{(k)}_i, Y^{(k)}_i} \big).
\eens
Now, on event $\MG$,
we have by Lemma~\ref{lemma::eventME},
simultaneously for all $\Delta'_{\Omega}$ as in \eqref{eq::diagsum},
\bens
\nonumber
\lefteqn{  \abs{\ip{\diag(\Delta'_{\Omega}), \hat{S} - \Sigma_0} } \le
\sum_{k} \abs{\ip{\diag( \tilde{\Delta}_k), \hat{S} - \Sigma_0}}} \\
\label{eq::aggre}
& \le & \sum_{k} \frac{t_{k} \twonorm{\delta^k}}{1-\ve} =
\sum_{k} C_m \twonorm{\Sigma_0} \fnorm{\diag(\Delta'_{\Psi_k})}
 (\sqrt{p} \vee d_k ).
 \eens
By the bound immediately above and \eqref{eq::diagsumproof},  we
obtain on event $\MG \cap \drp_0$,
\bens
\lefteqn{\abs{\ip{\diag(\Delta_{\Omega}), \hat{S} - \Sigma_0} }
\le  \abs{\ip{\tau_p     I_p, \hat{S} - \Sigma_0}} + 
\abs{\ip{\diag(\Delta'_{\Omega}), \hat{S} - \Sigma_0} } }\\
& \le & C_0  \twonorm{\Sigma_0}\big( \tau_{\Omega} \sqrt{p \log p} +  
  \sum_{k=1}^L \sqrt{p}  \fnorm{\diag(\Delta'_{\Psi_k})} \big)\\
&& +  C_m \twonorm{\Sigma_0} \sum_{k=1}^L d_k
\fnorm{\diag(\Delta'_{\Psi_k})} =:  r_{\diag, 1} + r_{\diag, 2},
\eens
where by Lemma~\ref{lemma::orthoDelta}, for ${r_{\diag,2}}/{(C_m \twonorm{\Sigma_0} )}$,
\bens
\sum_{k=1}^L d_{k} \fnorm{\diag(\Delta'_{\Psi_k})}
\le \sqrt{ \frac{d_{\max}}{m_{\min}}} \sqrt{L d_{\max}} 
\fnorm{\diag(\Delta_{\Omega})},
\eens
and by the Cauchy-Schwarz inequality,
\bens
\lefteqn{\frac{r_{\diag,1}}{C_0 \twonorm{\Sigma_0} }
  :=\tau_{\Omega} \sqrt{p} \sqrt{\log p} +\sum_{k=1}^L \sqrt{d_k}
  \sqrt{m_k}  \fnorm{\diag(\Delta'_{\Psi_k})} } \\
&  \le &   
  \big(\log p+ \sum_{k=1}^L d_k \big)^{1/2} \big(\sum_{k=1}^L m_k  \fnorm{\diag(\Delta'_{\Psi_k})}^2 + \tau_{\Omega}^2 p \big)^{1/2}  \\
  &  \le &   c  \big(\sum_{k=1}^L d_k \big)^{1/2}
  \fnorm{\diag(\Delta_\Omega)} \le  c \sqrt{L d_{\max}} \fnorm{\diag(\Delta_{\Omega})},
\eens
where $\log p = \sum_{k=1} \log d_k \le \sum_{k=1}^L d_k$, since the
RHS is a polynomial function of $p$, and the last line holds by
~\eqref{eq::ortho}.
Putting things together, we have 
\bens
\frac{r_{\diag}} {\twonorm{\Sigma_0}}
& \le &  
C_1 \sqrt{d_{\max}} \sqrt{L} \fnorm{\diag(\Delta_{\Omega})} 
\big(1 \vee \sqrt{{d_{\max}}/{m_{\min}}} \big).
\eens
To bound $\drp_0$,  we rewrite the trace as a quadratic form:
\bens 
\ip{\hat{S} - \Sigma_0, I} = \tr(\hat{S} - \Sigma_0) = Z^T \Sigma_0 Z
-\E (Z^T \Sigma_0 Z),
\eens 
where $Z \in \R^p$ is the same as in~\eqref{eq::tensordata}.
Thus, we have by the Hanson-Wright inequality \cite{RV13},
cf. Theorem 1.1 therein, and $\fnorm{\Sigma_{0} } \le \sqrt{p} \twonorm{\Sigma_0 }$,
\bens
\lefteqn{\prob{\abs{\ip{\hat{S} - \Sigma_0, I} > C  \twonorm{\Sigma_0 }
      \sqrt{p \log p}} }} \\
& \le & 
2 \exp\big(-c\min\big(\frac{C^2 p \log p \twonorm{\Sigma_0}^2}{\fnorm{\Sigma_0 }^2}, 
    C \sqrt{p \log p}\big) \big)  \le  \inv{p^4},
\eens
where $(C^2 \wedge C) c \ge 4$.
Hence by Lemma~\ref{lemma::eventME} and 
the bound immediately above, 
\bens 
\prob{\MG \cap \drp_0} \ge 1- c' \exp(-\log p) - \sum_{k} \exp(-c d_k).
\eens
The lemma thus holds upon adjusting the constants.
\end{proofof2}

\subsection{Proof of Lemma~\ref{lemma::eventME}}
\label{sec::eventME}
\begin{proofof2}
Set  $t_{k} > 0$.
First, we rewrite~\eqref{eq::YYtrace} and the trace term as a quadratic form in subgaussian random variables,
\ben
\label{eq::WZ}
&& \quad \ip{\diag( \tilde{\Delta}_k), \hat{S} - \Sigma_0}
=Z^T W Z - \E (Z^T W
Z), \\
\nonumber
&& \text{with $Z \in \R^p$ as in~\eqref{eq::tensordata} and }
W:= \Sigma_0^{1/2} \diag(\tilde{\Delta}_k) \Sigma_0^{1/2}.
\een
Then $\norm{W} \le \norm{\diag(\tilde{\Delta}_k)}\twonorm{\Sigma_0}$, 
where $\norm{\cdot}$ represents the  operator or the Frobenius norm.
Now  for $\delta \in \R^{d_k}$, by \eqref{eq::YYtrace},~\eqref{eq::WZ}, and the Hanson-Wright inequality,
\ben
\nonumber
\prob{\abs{\sum_{i=1}^{d_k}    \frac{\delta_{i}}{\twonorm{\delta}}
      \big(\twonorm{Y^{(k)}_i}^2 - \E \twonorm{Y^{(k)}_i}^2 \big)} \ge
    t_{k}}
  & = &
  \nonumber
\prob{\abs{\ip{\diag( \tilde{\Delta}_k), \hat{S} -
      \Sigma_0}}\ge t_{k} {\twonorm{\delta}}} \\
\nonumber
& = &
\prob{\abs{Z^T W
    Z - \E  (Z^T W Z)} \ge t_{k} {\twonorm{\delta}}}    \\
\nonumber
& \le &  
2 \exp\left[-c\min\big(\frac{t_{k}^{2} \twonorm{\delta}^2}{
      \fnorm{W}^2},     \frac{t_{k} {\twonorm{\delta}}}{\twonorm{W}}
    \big) \right] \\
  \label{eq::crossB}
  & =: & p_1.
\een
Now for all $\delta = (\delta_{1},  \ldots, \delta_{d_k})$ and
$\diag(\Delta'_{\Psi_k}) = \diag(\delta)$, we have
\bens
\twonorm{\diag(\tilde{\Delta}_k)} & = &
\twonorm{\diag(\Delta'_{\Psi_k})} \; \text{  and} \\
\fnorm{\diag(\tilde{\Delta}_k)} & = &
\sqrt{m_k}
\fnorm{\diag(\Delta'_{\Psi_k})} = \sqrt{m_k} \twonorm{\delta}
\eens
by~\eqref{eq::diagDK}.
Thus
\bens
&& \twonorm{W} \le \twonorm{\Sigma_0} \twonorm{\diag(\tilde{\Delta}_k)}
\le  \twonorm{\Sigma_0} \twonorm{\delta}, \; \text{and  } \\
&& 
\fnorm{W} \le \twonorm{\Sigma_0} \fnorm{\diag(\Delta'_{\Psi_k})} 
= {\twonorm{\Sigma_0} }\sqrt{m_k} \twonorm{\delta}.
\eens
Recall $\Pi_{d_k}$ is an $\ve$-net of the sphere $\Sp^{d_k-1}$, where
$0< \ve < 1/2$.
Then for $t_{k}   := C_m
\twonorm{\Sigma_0}  (\sqrt{p} \vee d_k)$ as in \eqref{eq::tk},
we have by~\eqref{eq::crossB} and the union bound,
\bens
\lefteqn{\prob{\exists \delta \in   \Pi_{d_k}: \sum_{i=1}^{d_k}
    \delta_{i}    \big(\twonorm{Y^{(k)}_i}^2 - \E \twonorm{Y^{(k)}_i}^2 \big) \ge
    t_{k}} }\\
&=:&  \prob{\text{ event } \; \MG^c_k \; \text{occurs } } \le
(1+ 2/\ve )^{d_k} p_1 \\
&\le &
5^{d_k}  \exp\big(-c\min\big(\frac{ C_m^2
      \twonorm{\Sigma_0}^2 p}{m_k  \twonorm{\Sigma_0}^2},
    \frac{ C_m\twonorm{\Sigma_0}  d_k}{ \twonorm{\Sigma_0}}\big) \big) \\
  & \le & \exp(d_k \log 5 -c d_k (C_m^2 \wedge C_m))\le   \exp(-c' d_k \log 5).
\eens
The ``moreover'' statement follows from a standard
approximation argument. Suppose event $\MG$ holds.
\silent{Then on event $\MG_k$, for a chosen parameter $t_{k}$ as 
in  \eqref{eq::tk}, we have
\bens
\lefteqn{
  \sup_{\delta \in \Sp^{d_k-1}}
  \sum_{i=1}^{d_k} \delta_{i}
  \big(\ip{Y^{(k)}_i, Y^{(k)}_i} - \E \ip{Y^{(k)}_i, Y^{(k)}_i} \big)} \\
& \le &
\inv{1-\ve} \sup_{\delta \in \Pi_{d_k}} \sum_{j=1}^{d_k}
\delta_{i}  \big(\ip{Y^{(k)}_i, Y^{(k)}_i} - \E \ip{Y^{(k)}_i,
    Y^{(k)}_i}\big) \\
& \le & \frac{t_{k}}{1-\ve},
\eens
by a standard approximation argument and the union bound.}
Denote by
\bens
y =\big(\twonorm{Y_1^{(k)}}^2 - \E\twonorm{Y_1^{(k)}}^2,
  \ldots,   \twonorm{Y_{d_k}^{(k)}}^2 - \E\twonorm{Y_{d_k}^{(k)}}^2\big).
\eens
We have for $\delta = (\delta_1, \ldots, \delta_{d_k}) \in \Sp^{d_k-1}$,
\bens
\sup_{\delta \in \Pi_{d_k}} \ip{\delta, y}
\le \twonorm{y} = \sup_{\delta \in \Sp^{d_k-1}} \ip{y, \delta}
\le \inv{1-\ve} \sup_{\delta \in \Pi_{d_k}} \ip{\delta, y}.
\eens
The LHS is obvious.
To see the RHS, notice that for $\delta \in \Sp^{d_k-1}$ that achieves maximality in
\bens
\twonorm{y} = \sup_{\delta \in \Sp^{d_k-1}} \ip{y, \delta},
\eens
we can find $\delta_0 \in \Pi_{d_k}$ such that $\twonorm{\delta -
  \delta_0} \le \ve$.
Now
\bens
\ip{\delta_0, y} & = &  \ip{\delta, y} - \ip{\delta-\delta_0, y} \\
  & \ge &  \ip{\delta, y} - \sup_{\delta \in \Sp^{d_k-1}} \ve
  \ip{\delta, y} =  (1-\ve) \sup_{\delta \in \Sp^{d_k-1}} \ip{\delta,   y},
\eens
\text{and hence}
\bens
\sup_{\delta \in \Pi_{d_k} } \ip{\delta, y}
  & \ge &  (1-\ve) \sup_{\delta
    \in \Sp^{d_k-1}} \ip{\delta, y} =  (1-\ve) \twonorm{y}.
\eens
The lemma thus holds.
\end{proofof2}

\section{Proof of Theorem~\ref{thm::main}}
\label{sec::strategy}

First we state Theorem \ref{thm::orig} from~\cite{GZH19}.
  
\begin{theorem}[\cite{GZH19}, restated]
  \label{thm::orig}
Suppose (A1) and (A2) hold and $n(m_{\min})^2 \ge C^2
\kappa(\Sigma_0)^4 (s+ p) (L+1)^2 \log p$, where $s = \sum_{k} m_k s_k$  is as in Definition~\ref{def::A1}.
Then
\bens
\label{eq::oldFrate}
\frac{\|\hat{{\Omega}} - {\Omega}_0\|_F}{\|\Omega_0\|_2} &= &
O_p\big(  \kappa(\Sigma_0) \sqrt{L+1}
  \big(\frac{(s + p)\log p}{n m_{\min}}\big)^{1/2}\big),\\
\label{eq::oldrate}
\frac{\|\hat{{\Omega}} - {\Omega}_0\|_2}{\|\Omega_0\|_2} &= & O_p
\big( \kappa(\Sigma_0) (L+1) \big(\frac{(s +  p)\log p}{n m^2_{\min}}\big)^{1/2}\big).
\eens
\end{theorem}

Recall \eqref{eq::objfunc} is equivalent to
\bens
\label{eq::lossfunc}
&& \hat{\Omega}
= \arg\min_{\Omega \in \mathcal{K}_{\mathbf{p}}^\sharp}
\big(-\log \abs{\Omega} + \ip{\hat{S}, \Omega} +
  \sum_{k=1}^L  m_k \rho_{n,k} \abs{{\Psi}_k}_{1,\off}\big),
\eens
where $\hat{S}$ is as defined in \eqref{eq::gram}, in view of
\eqref{eq::gramK2}.
First, we define the unified event $\mathcal{A}$ as the event that all these events hold, i.e.
\bens
  \mathcal{A} = \T \cap \mathcal{D}_0 \cap \MG,
\quad  \text{ where} \quad \mathcal{G} = \mathcal{G}_1 \cap \dots \cap \mathcal{G}_L.
\eens
We focus on the case $n=1$. For $n >1$, we defer the proof to
Section~\ref{sec::mainthmlargen}.
First, we state Lemma~\ref{thm::pivot}, which is proved in
~\cite{GZH19supp}, cf. Lemma 8 therein.
\begin{lemma}{\textnormal{(Lemma 8 of~\cite{GZH19supp})}}
\label{thm::pivot}
  For all $\Omega \in \mathcal{K}_{\bp}$, 
$\twonorm{\Omega} \le \sqrt{\frac{L+1}{ \min_k m_k}}\fnorm{\Omega}$. 
\end{lemma}

In the proof of Theorem
\ref{thm::main} that follows, our strategy will be to show that
several events controlling the concentration of the sample covariance
matrix (in the $n=1$ case, simply an outer product) hold with high
probability, and then show that given these events hold, the
statistical error bounds in Theorem \ref{thm::main} hold.
The off-diagonal events are as defined in~\eqref{eq::defineToffd}.

We adopt the definitions of new diagonal events in Section~\ref{sec::proofofdiagfinal}.
We use the following notation to describe errors in the precision matrix and its factors. For $\Omega \in \ksump$ let $\Delta_{\Omega} = \Omega - \Omega_0 \in \ksump$. Since both $\Omega$ and $\Omega_0$ are Kronecker sums,
\bens 
\Delta_\Omega &  = &  \Delta_{\Psi_1}
\oplus \Delta_{\Psi_2}  \oplus \ldots \oplus \Delta_{\Psi_L}
\eens
for some $\Delta_{\Psi_k}$ whose off-diagonal (but not diagonal) elements are uniquely determined. For an index set $S$ and a matrix $W = [w_{ij}]$, write $W_S \equiv
(w_{ij} I( (i,j) \in S))$, where $I(\cdot)$ is an indicator function.

\subsection{Preliminary results}
Before we show the proof of Theorem~\ref{thm::main},
we need to state the following lemmas.
We then present an error bound for the off-diagonal component of the loss function,
which appears as Lemma 12 in~\cite{GZH19supp} and follows from the
concentration of measure bounds on elements of $\offd(S^k -
\Sigma_0^{(k)})$; cf.~\eqref{eq::defineToffd}.
Combined with our new concentration bound on the diagonal component of
the loss function, cf. Lemma~\ref{lemma::diagnew}, we obtain the improved overall rate of convergence as stated in Theorem~\ref{thm::main}. 
\begin{lemma}
\label{lemma::triangle}
Let $\Omega_0 \succ 0$.
Let $S= \{ (i,j) : \ \Omega_{0ij} \neq 0, \ i\neq j \}$ and 
$S^c = \{ (i,j) : \ \Omega_{0ij} = 0, \ i\neq j \}$.
Then for all $\Delta \in \ksump$, we have
\begin{eqnarray} 
\label{eq::Delta-I}
\offone{\Omega_0 + \Delta}  - \offone{\Omega_0}
& \geq & \onenorm{\Delta_{S^c}} - \onenorm{\Delta_S}
\end{eqnarray}
where by disjointness of $\supp(\offd(\Psi_k)) := \{(i,j): i\not=j, \; \Psi_{k, ij} \not= 0\}, k=1, \ldots, L$,
\bens
\onenorm{\Delta_{S}}  = \sum_{k=1}^L m_k
\onenorm{\Delta_{\Psi_k,S}} \; \text{ and } \; 
\onenorm{\Delta_{\Sc}}  = \sum_{k=1}^L m_k
\onenorm{\Delta_{\Psi_k, \Sc}}.
\eens
\end{lemma}

Proofs of Lemmas \ref{thm::pivot} and \ref{lemma::triangle} appear 
in~\cite{GZH19supp} (cf.  Lemmas 8 and 10 therein). 
Lemma~\ref{lemma::offdold} follows from~\cite{GZH19supp}; cf. Lemmas 11 and~12 therein.
\begin{lemma}
\label{lemma::offdold}
With probability at least $1 - 2 L \exp(-c' \log p)$,
\bens
\label{eq::offdbounds}
&& \abs{\ip{\offd(\Delta_{\Omega}), \hat{S} - \Sigma_0} }\le \sum_{k=1}^L
m_k \offone{\Delta_{\Psi_k}} \delta_{k}, \\
&& \; \text{ where} \;
\delta_{k} \asymp \sqrt{\frac{\log p}{m_k}} \twonorm{\Sigma_0}, \forall k.
\eens
\end{lemma}

Next we show that as an immediate corollary of \eqref{eq::Delta-I},
we have Lemma~\ref{lemma::deltagbound}, which is a deterministic
result and identical to Lemma~10~\cite{GZH19supp}. The proof is
omitted.
\begin{lemma}{\textnormal{(Deterministic bounds)}}
  \label{lemma::deltagbound}
  Let $\rho_{k} \ge 0$.  Denote by
  \ben
  \label{eq::deltag}
  \Delta_g
  &  := &
  \sum_{k=1}^L m_k \rho_{k}
  \left(\offone{\Psi_k + \Delta_{\Psi_k}} - \offone{\Psi_k}\right), \\
  \nonumber
  \text{ then} \quad
\Delta_g  &  \ge &  \sum_{k=1}^L m_k \rho_{k}
\left(\onenorm{\Delta_{\Psi_k, \Sc}} - \onenorm{ \Delta_{\Psi_k,  S}}\right).
\een
\end{lemma}

Lemma~\ref{lemma::offdpre} follows immediately from 
Lemmas~\ref{lemma::offdold} and~\ref{lemma::deltagbound}. 
\begin{lemma}
  \label{lemma::offdpre}
  Suppose that $d_k = O(m_k)$ for all $k$. Let  $\Delta_g$ be as in Lemma~\ref{lemma::deltagbound}.
  Under the settings of Lemmas~\ref{lemma::offdold} and~\ref{lemma::deltagbound},
  we have for choices of $\rho_{k} = \delta_{k}/\ve_k, \forall k$, where $0 <
  \ve_k < 1$ and $\delta_{k} \asymp \sqrt{\frac{\log p}{m_k}}  \twonorm{\Sigma_0}$,
\ben
\label{eq::offdbound}
\Delta_g +\ip{\offd(\Delta_{\Omega}), \hat{S} - \Sigma_0} \ge
-2 \max_{k} \rho_{k} \onenorm{\Delta_{S}}.
 \een
\end{lemma}

\begin{proof}
First, we prove~\eqref{eq::offdbound}.
We have by \eqref{eq::deltag}
\bens
\lefteqn{\Delta_g +\ip{\offd(\Delta_{\Omega}), \hat{S} - \Sigma_0} } \\
& \ge  &  \sum_{k=1}^L m_k \rho_{k}
\left(\offone{\Psi_k + \Delta_{\Psi_k}} - \offone{\Psi_k}\right) + 
\ip{\offd(\Delta_{\Omega}), \hat{S} - \Sigma_0} =: S_2
\eens
where under the settings of Lemma~\ref{lemma::offdold},
\bens
S_2
 & \ge & 
\sum_{k=1}^L m_k \rho_{k}
\left(\onenorm{\Delta_{\Psi_k, \Sc}} - \onenorm{ \Delta_{\Psi_k,
      S}}\right) -  \sum_{k=1}^L m_k 
\offone{\Delta_{\Psi_k}}
\delta_{k} \\
\nonumber
 & \ge & 
\sum_{k=1}^L m_k \rho_{k}
\left(\onenorm{\Delta_{\Psi_k, \Sc}} - \onenorm{ \Delta_{\Psi_k,
      S}}\right)  -  \sum_{k=1}^L m_k \delta_{k} 
\left(\onenorm{\Delta_{\Psi_k, \Sc}} + \onenorm{\Delta_{\Psi_k,
      S}}\right) \\
 & \ge & 
 -\sum_{k=1}^L m_k (\rho_{k} + \delta_{k})
 \onenorm{\Delta_{\Psi_k, 
     S}}  \\
& \ge & 
-2 \max_{k} \rho_{k} \sum_{k=1}^L m_k \onenorm{\Delta_{\Psi_k,  S}} = 
-2 \max_{k} \rho_{k} \onenorm{\Delta_{S}};
\eens
Thus \eqref{eq::offdbound} holds.
\end{proof}

Lemma~\ref{lemma::threeamigo}
follows from Lemmas~\ref{lemma::diagnew}
and~\ref{lemma::offdpre}.  
We defer the proof of
Lemma~\ref{lemma::threeamigo} to Section~\ref{sec::proofofT3}.
Since $p = \prod_{k} d_k \ge 2^{L}$ so long as $d_k \ge 2$, we have $\log p \ge L$ and hence $\exp(c \log p) > L$ for sufficiently large $c$.
\begin{lemma}
    \label{lemma::threeamigo}
    Suppose that $n=1$.  Let $s = \sum_{k=1}^L m_k s_k$.
    Then, under the settings of Lemmas~\ref{lemma::diagnew} 
and~\ref{lemma::offdpre}, we have with probability at least $1 - L \exp(-c' \log p)$,
\bens
\abs{\Delta_g + \ip{\Delta_{\Omega}, \hat{S} - \Sigma_0} }
& \le &
C' \twonorm{\Sigma_0} T_3 \; \text{ where} \; 
T_3  := \frac{\sqrt{s \log p + L p} \fnorm{\Delta_{\Omega}}}{\sqrt{m_{\min}}}.
\eens
\end{lemma}

\begin{proposition}
  \label{prop::prelude}
  Set $C > 36  (\max_{k} \inv{\ve_k} \vee 
C_{\diag})$ for $C_{\diag}$ as in  Lemma~\ref{lemma::diagfinal}. 
 Let
  \ben 
\label{eq:Mdef}
r_{\mathbf{p}} &  = & C \twonorm{\Sigma_0} \sqrt{s \log p + Lp
}/\sqrt{m_{\min}}  \text{ and } \quad M   = \frac{1}{2}{\phi^2_{\max}(\Omega_0)}
= \inv{2 \phi^2_{\min}(\Sigma_0)}.
\een
Let  $\Delta_\Omega \in \mathcal{K}_{\bp}$ such that
$\fnorm{\Delta_\Omega} = M r_{\mathbf{p}}$.
Then $\twonorm{\Delta_{\Omega}} \le \half \phi_{\min}(\Omega_0).$
\end{proposition}

\begin{proof}
Indeed, by Theorem~\ref{thm::pivot}, we have for all $\Delta \in \mathcal{T}_n$,
\bens
\twonorm{\Delta}  & \leq & 
\sqrt{\frac{L+1}{\min_k m_k}} \|\Delta\|_F  =
\sqrt{\frac{L+1}{m_{\min}}} M r_{\mathbf{p}}  \\
& \le & 
\sqrt{\frac{L+1}{m_{\min}}} \frac{C}{2} \inv{\phi_{\min}^2(\Sigma_0)}
\twonorm{\Sigma_0} \sqrt{\frac{s  \log p + p  L}{m_{\min}}}  \le  \half \phi_{\min}(\Omega_0)  = \inv{2 \phi_{\max}(\Sigma_0)} 
\eens
so long as $m_{\min}^2 > 2 C^2 (L+1)\kappa(\Sigma_0)^4 (s \log p+ p L)$,
where $\kappa(\Sigma_0)$ is the condition number of $\Sigma_0$.
\end{proof}

\subsection{Proof of Theorem~\ref{thm::main}}
\begin{proofof2}
  We will only show the proof for $n=1$.
Let
\ben
  \label{eq:Qbar}
G(\Delta_\Omega)& =&Q(\Omega_0 + \Delta_\Omega) - Q(\Omega_0)
 \een                  
be the difference between the objective function~\eqref{eq::lossfunc} at $\Omega_0 + \Delta_\Omega$ and at $\Omega_0$. Clearly $\hat{\Delta}_\Omega = \hat{\Omega}-\Omega_0$ minimizes
$G(\Delta_\Omega)$, which is a convex function  with a unique
minimizer on $\mathcal{K}_{\mathbf{p}}^\sharp$ (cf. Theorem
5~\cite{GZH19supp}).  Let $r_{\bp}$ be as defined in \eqref{eq:Mdef} for some large enough
  absolute constant $C$ to be specified, and
\begin{equation}\label{eq:rrnM}
\T_n =  \left\{\Delta_\Omega \in \mathcal{K}_{\mathbf{p}}:
  \Delta_\Omega = \Omega - \Omega_0, \Omega,\Omega_0 \in
  \mathcal{K}_{\mathbf{p}}^\sharp, \|\Delta_\Omega\|_F = M
  r_{\mathbf{p}}\right\}.
\end{equation}
In particular, we set $C > 36  (\max_{k} \inv{\ve_k} \vee C_{\diag})$
in $ r_{\mathbf{p}}$, for absolute constant $C_{\diag}$ as in  Lemma~\ref{lemma::diagfinal}.
Proposition~\ref{prop:cnv} follows from~\cite{ZLW08}.
\begin{proposition}\label{prop:cnv}
If $G(\Delta) > 0$ for all $\Delta \in \mathcal{T}_n$ as defined in \eqref{eq:rrnM}, then $G(\Delta) > 0$ for all $\Delta$ in
\bens
\mathcal{V}_n = \{\Delta \in \mathcal{K}_{\mathbf{p}}: \Delta = \Omega - \Omega_0, \Omega,\Omega_0 \in \mathcal{K}_{\mathbf{p}}^{\sharp}, \|\Delta\|_F > M r_{\mathbf{p}}\}
\eens
for  $r_{\mathbf{p}}$~\eqref{eq:Mdef}.
Hence if $G(\Delta) > 0$ for all $\Delta \in \T_n$, then $G(\Delta) > 0$ for all $\Delta \in \T_n \cup \mathcal{V}_n$.
\end{proposition}

\begin{proposition}
\label{prop:bnd}
Suppose $G(\Delta_\Omega) > 0$ for all $\Delta_\Omega \in
\mathcal{T}_n$.
We then have
\bens
\fnorm{\hat{\Delta}_{\Omega}} < M r_{\mathbf{p}}.
\eens
\end{proposition}

\begin{proof}
By definition, $G(0) = 0$, so $G(\hat{\Delta}_\Omega) \leq G(0) =
0$. Thus if $G(\Delta_\Omega) > 0$ on $\mathcal{T}_n$, then by
Proposition \ref{prop:cnv}, $\hat{\Delta}_\Omega \notin \mathcal{T}_n \cup \mathcal{V}_n$ where
$\mathcal{V}_n$ is defined therein. The proposition thus holds.
\end{proof}

\begin{lemma}
\label{cor:LD} 
Under (A1) -  (A3), for all $\Delta \in {\mathcal{T}_n}$ for
which $r_{\mathbf{p}} = o\left(\sqrt{\frac{\min_k
      m_k}{L+1}}\right)$,
\bens
\log|\Omega_0 + \Delta| - \log|\Omega_0| \leq \langle \Sigma_0, \Delta
\rangle -
\frac{2}{9\|\Omega_0\|_2^2}\fnorm{\Delta}^2.
\eens
\end{lemma}
We defer the proof of Lemma~\ref{cor:LD} to Section~\ref{app:LD}.
By Proposition \ref{prop:bnd}, it remains to show that $G(\Delta_\Omega) > 0$ on
$\mathcal{T}_n$ under the settings of Lemma \ref{lemma::threeamigo}.
\begin{lemma}
  With probability at least $1 - L \exp(-c' \log p)$,
  we have $G(\Delta) > 0$ for all $\Delta \in \mathcal{T}_n$.
\end{lemma}

\begin{proof}
By Lemma \ref{cor:LD}, if $r_{\mathbf{p}} \leq \sqrt{\min_k m_k/(L+1)}$, 
we can express \eqref{eq:Qbar} as
\ben
\nonumber
\lefteqn{G(\Delta_\Omega) =
\langle\Omega_0 + \Delta_\Omega,\hat{S}\rangle - \log|\Omega_0 +
\Delta_\Omega| - \langle\Omega_0,\hat{S}\rangle + \log|\Omega_0|  }\\
\nonumber
& & + \underbrace{\sum_k \rho_k m_k (|\Psi_{k,0}+ \Delta_{\Psi,k}|_{1
    ,\off}- |\Psi_{k,0}|_{1,\off})}_{\Delta_g}  \\
& \geq&  
\label{eq:G}
\quad \quad \quad
\ip{\Delta_\Omega,\hat{S} -\Sigma_0}+ \frac{2}{9\|\Omega_0\|_2^2}
\|\Delta_\Omega\|_F^2 +  \Delta_g.
\een
By Lemma \ref{lemma::threeamigo} and \eqref{eq:G},
we have  for all $\Delta_\Omega \in \T_n$, and
$C' =\max_{k}(\frac{2}{\ve_k}) \vee 2 C_{\diag}$, 
\bens
\label{eq:G3}
G(\Delta_\Omega) 
&\geq & \frac{2}{9\|\Omega_0\|_2^2} \|\Delta_\Omega\|_F^2
-\abs{{\Delta_g} + \ip{\Delta_\Omega,\hat{S}- \Sigma_0 }} \\
&\geq & 
\frac{2}{9\twonorm{\Omega_0}^2} \fnorm{\Delta_\Omega}^2
- \frac{C' \twonorm{\Sigma_0}}{\sqrt{\min_k m_k }}
\sqrt{s\log p + Lp} \fnorm{\Delta_{\Omega}} =: W,
 \eens
 where by Lemma~\ref{lemma::threeamigo}, we have with probability at
 least $1 - L \exp(-c' \log p)$, 
\bens
\abs{\Delta_g + \ip{\Delta_{\Omega}, \hat{S} - \Sigma_0} }
& \le &
C' \twonorm{\Sigma_0}  \frac{\sqrt{s \log p + L p} \fnorm{\Delta_{\Omega}}}{\sqrt{m_{\min}}}
\eens
for $d_k =O(m_k)$.
Now $W >0$ for $\fnorm{\Delta_{\Omega}} = M r_{\mathbf{p}}$, where
$M = \inv{2 \phi_{\min}^2(\Sigma_0)}$, since
\bens
&& C' \twonorm{\Sigma_0} \sqrt{\frac{1}{\min_k m_k }} \sqrt{(L p+ s \log
  p)} \inv{M r_{\mathbf{p}}}
=   \frac{C'}{ C M} \\
&& =   \frac{2 C'}{C} \phi_{\min}^2(\Sigma_0) < \frac{2}{9\twonorm{\Omega_0}^2},
\eens
which holds so long as
$C$ is  chosen to be large enough in $r_{\mathbf{p}}$ as defined in~\eqref{eq:Mdef}.
For example, we set $C = 18 C' = 36 (\max_{k}(\frac{1}{\ve_k}) \vee  C_{\diag})$.
\end{proof}

Theorem~\ref{thm::main} follows from Proposition \ref{prop:bnd} immediately.
Combining Lemmas \ref{lemma::offdold} and \ref{lemma::diagnew} using
the union bound implies both events hold with probability at least $ 1
- L \exp(-c' \log p)$.
The error in the operator norm immediately follows from the Frobenius norm error
bound and Lemma~\ref{thm::pivot}.
\end{proofof2}

To complete the proof, it remains to present the case of $n >
1$. We leave the details to Section~\ref{sec::mainthmlargen} for completeness.

\section{Proof of preliminary results in Section~\ref{sec::strategy}}
\label{sec::proofofmainsupp}


\subsection{Proof of Lemma~\ref{lemma::threeamigo}}
\label{sec::proofofT3}
\begin{proofof2}
We focus on the case $d_k \le m_k \forall k$;
By definition of $\Delta_g$,
\bens
\ip{\Delta, S - \Sigma_0} + \Delta_g
& := &  \ip{\offd(\Delta), S - \Sigma_0} + \Delta_g + \\
&& 
 \ip{\diag(\Delta), S -   \Sigma_0} 
 \eens
Then we have by \eqref{eq::offdbound} and \eqref{eq::diagfinal},
with probability at least
$$1-\sum_{k=1}^L 2 \exp(-c d_k)  - 2 L \exp(-c' \log p), \; \text{ for
} \; d_k = O(m_k),$$
and $\onenorm{\Delta_{S}} \le \sqrt{s} \fnorm{\Delta_{S}},$
where $s = \sum_{k=1}^L m_k s_k$, 
\bens
\abs{\Delta_g + \ip{\offd(\Delta_{\Omega}), \hat{S} - \Sigma_0} }
& \le &
2 \max_{k} \rho_{k} \onenorm{\Delta_{S}} \\
& \le &
 2 \max_{k} \big(\inv{\ve_{k}} \sqrt{\frac{\log p}{m_k}}\big)
 \sqrt{s} \fnorm{\Delta_{\Omega, S}}
 \eens
and
\bens
\abs{\ip{\diag(\Delta_\Omega), \hat{S}-\Sigma_0} }
& \leq &  C_{\diag} \twonorm{\Sigma_0} \sqrt{d_{\max} } \sqrt{L} \fnorm{\diag(\Delta_{\Omega})} \big(1 + \sqrt{\frac{d_{\max}}{m_{\min}}} \big).
\eens
Let $C_{\offd} := \max_{k} \big(1/{\ve_k} \big)$ and $C' = 2
(C_{\diag} \vee  C_{\offd})$, where $C_{\offd} = 2 \max_{k} \inv{\ve_k}$.
The Lemma thus holds by the triangle inequality: for $d_{\max} \le \sqrt{p}$
\bens 
\lefteqn{
  \abs{\ip{\Delta, S - \Sigma_0} + \Delta_g}
  \leq \abs{ \ip{\offd(\Delta), \hat{S} - \Sigma_0} + \Delta_g}  +
 \abs{\ip{\diag(\Delta), \hat{S} -  \Sigma_0} }}\\
 & \le &
2 C_{\offd} \twonorm{\Sigma_0} \sqrt{\frac{s \log p}{m_{\min}}}
\fnorm{\Delta_{\Omega, S}} + \\
&&
C_{\diag} \twonorm{\Sigma_0} \sqrt{d_{\max} } \sqrt{L} \fnorm{\diag(\Delta_{\Omega})} 
\left(1 + \sqrt{\frac{d_{\max}}{m_{\min}}} \right) \\
 & \le &
2 C_{\offd} \vee C_{\diag} \twonorm{\Sigma_0} \big(\sqrt{\frac{ s \log p}{m_{\min}}} 
\fnorm{\Delta_{\Omega,S}} + \sqrt{L} \fnorm{\diag(\Delta_{\Omega})} 
\frac{\sqrt{p} + d_{\max}}{2 \sqrt{m_{\min}}} \big)\\
& \le &
C' \twonorm{\Sigma_0} T_3  
\eens  
where by Cauchy-Schwarz inequality,
\bens
\sqrt{s \log p} \fnorm{\offd(\Delta_{\Omega})}
+ \sqrt{L p} \fnorm{\diag(\Delta_{\Omega})} \le 
\sqrt{s \log p + p L} \fnorm{\Delta_{\Omega}}.
\eens
\end{proofof2}

\subsection{Proof of Lemma \ref{cor:LD}}
\label{app:LD}
\begin{proofof2}
We first state Proposition \ref{prop:posi-def-interval}
\begin{proposition}
\label{prop:posi-def-interval}
Under (A1)-(A3), for all $\Delta \in {\mathcal{T}_n}$,
\ben
\label{eq::eigen-bound-local}
\twonorm{\Delta} \le  M r_{\mathbf{p}}
\sqrt{\frac{L+1}{\min_k  m_k}} \le \half \phi_{\min} (\Omega_0),
\een
so that $\Omega_0 + v \Delta \succ 0, \forall v \in I \supset [0, 1]$,
where $I$ is an open interval containing $[0, 1]$.
\end{proposition}
\begin{proof}
  By Proposition~\ref{prop::prelude},~\eqref{eq::eigen-bound-local} holds for $\Delta \in {\mathcal{T}_n}$;
Next, it is sufficient to show that $\Omega_0 + (1 + \ve) \Delta \succ 0$
and $\Omega_0 - \ve \Delta \succ 0$ for some $1 > \ve > 0$.
Indeed,  for $\ve < 1$, 
\bens
\phi_{\min} (\Omega_0 + (1 + \ve) \Delta) 
& \geq &
\phi_{\min} (\Omega_0) - (1 + \ve) \twonorm{\Delta} \\
& > & \phi_{\min} (\Omega_0) - 2 \sqrt{\frac{L+1}{\min_k m_k}}M
r_{\mathbf{p}}   >0
\eens
given that by definition of $\T_n$ and \eqref{eq::eigen-bound-local}.
\end{proof}

Thus we have that $\log|\Omega_0 + v \Delta|$ is infinitely differentiable on
the open interval $I \supset [0, 1]$ of $v$. This allows us to 
use the Taylor's formula with integral remainder to prove Lemma
\ref{cor:LD}, following identical steps in~\cite{GZH19supp},
drawn from~\cite{RBLZ08}, and hence is omitted.
\end{proofof2}

\subsection{Extension to multiple samples $n > 1$}
\label{sec::mainthmlargen}
Incorporating $n > 1$ directly into the proof above is relatively
straightforward but notation-dense; hence it suffices to note that having $n$ independent samples essentially increases the $m_k$ replication to $n m_k$, and propagate this fact through the proof. 
We also note that the multi-sample $n > 1$ case can be converted to
the single sample $n=1$ regime to obtain a result directly. To see
this, note that $n$ independent samples with precision matrix $\Omega_0
\in \mathbb{R}^{p \times p}$ can be represented as a single sample
with the block-diagonal precision matrix, i.e. $\Omega_0$ repeated $n$
times blockwise along the diagonal, specifically, $\Omega^{(n)} = I_n
\otimes \Omega_0 \in \mathbb{R}^{pn \times pn}$.
Recall that by definition of the Kronecker sum,
$$\Omega^{(n)} = I_n \otimes \Omega_0 =
0_{n\times n} \oplus \Psi_1 \oplus\dots\oplus \Psi_L$$
is a $(L+1)$-order Kronecker sum with $p^{(n)} = pn$,
achieved by introducing an all-zero factor $\Psi_0 = 0_{n\times n}$
with $d_0 = n$ (and $m_0 = p$).
Since this extra factor is zero, the operator norms are not affected.
The sparsity factor of $\Omega^{(n)}$ is $s^{(n)} = s n$ since the
non-zero elements are replicated $n$ times, and each co-dimension
$m^{(n)}_k := p^{(n)}/d_k = n m_k$ for $k > 0$.

Hence the single sample convergence result can be applied with
$L^{(n)} = L+1$, yielding for $n \leq d_{\max}$ and $L \geq 2$
\bens
\twonorm{\hat{\Omega} -  \Omega_0}/  \twonorm{\Omega_0
}  &= &
C \kappa(\Sigma_0)\sqrt{L^{(n)} + 1} \sqrt{\frac{s^{(n)} \log p^{(n)} + L^{(n)} p^{(n)}}{[m^{(n)}_{\min}]^2}}\\
&= & C \kappa(\Sigma_0)\sqrt{L + 2}  \sqrt{\frac{s (\log p + \log n) + (L + 1) p}{n m_{\min}^2}}\\
&\leq & C \sqrt{\frac{8}{3}} \kappa(\Sigma_0)\sqrt{L + 1} \sqrt{\frac{s \log p  + L p}{n m_{\min}^2}}
\eens
since $m^{(n)}_{\min} = \min(m_0, n m_{\min}) = \min(d_{\max} m_{\min}, n m_{\min})= n m_{\min}$ whenever $n \leq d_{\max}$.
Hence Theorem \ref{thm::main} is recovered for $n \leq d_{\max}$, with constant slightly worse than could be obtained by incorporating $n$ directly into the proof.

\subsection{Proof of Corollary \ref{thm::cubic}}
\label{sec::proofofcubic}
 \begin{proofof2}
   Denote by $s = \sum_{k} m_k s_k$.
   Then for $n=1$ and $p = d_{\max}
  m_{\min} = m_k d_k$ for all $k$,
   \bens
\sqrt{L}  \sqrt{\frac{s \log p + L p}{m_{\min}^2}}
   & = &
   \sqrt{L \frac{d_{\max}}{m_{\min}}}
   \sqrt{\frac{\sum_{k} m_k s_k \log p +      L p}{ d_{\max} m_{\min}}} \\
& = &
   \sqrt{L \frac{d_{\max}}{m_{\min}}}   \sqrt{\frac{\sum_{k} m_k s_k
       \log p +    L p}{p}}\\
   & = & \sqrt{L \frac{d_{\max}}{m_{\min}}}
   \sqrt{\sum_{k} \frac{s_k\log p}{d_k} + L} <1
\eens
by (A3); The corollary thus follows from Theorem~\ref{thm::main}.
\end{proofof2}

\subsection{Proof of Theorem~\ref{thm::maincubic}}
\begin{proofof2}
  \label{sec::proofmaincubic}
  Suppose that $m_1\asymp m_2 \asymp \ldots \asymp m_L$.
   Denote by $s = \sum_{k} m_k s_k$. Then
  \bens
 \sqrt{\frac{s \log p + L p}{(\min_{k} m_k)}}
  & = &
  \sqrt{\frac{\sum_{k} m_k s_k \log p + L p}{ m_{\min}}} \\
&   \approx &
\sqrt{L}    \sqrt{\inv{L} \sum_{k} s_k \log p + d_{\max}}
\eens
The theorem thus follows from Theorem~\ref{thm::main}.
\end{proofof2}

\section{Conclusion}
\label{sec::conclude}
We present sharper statistical rates of
convergence of the $\ell_1$ regularized TeraLasso estimator of
precision matrices with Kronecker sum structures in the finite sample settings.
The key innovation in the present work
is to derive tight concentration bounds for the trace terms on the diagonal
component of the loss function \eqref{eq::objfunc}.
Crucially, this improvement allows for finite sample
statistical rates of convergence to be derived for the two-way Kronecker sum model, which was
missing from~\cite{GZH19} and was also deemed as the most demanding, due to the
lack of {\em sample replications} in complex and high-dimensional data.

\section*{Acknowledgement}

We thank Harrison Zhou for helpful discussions.  We thank the Simons Institute
for the Theory of Computing at Berkeley on the
occasion of Algorithmic Advances for Statistical Inference with Combinatorial
Structure Workshop, and organizers of International Conference on Statistics
and Related Fields (ICON STARF), University of Luxembourg, for their kind
invitations, where we presented a talk including this work in 2021.

\bibliography{subgaussian,KronML_bib}

\end{document}